\documentclass[12pt]{amsart}
\usepackage{hyperref}
\usepackage[utf8]{inputenc}
\usepackage{amssymb, graphicx}
\usepackage{amsfonts}
\usepackage{amsmath}
\usepackage{latexsym}
\usepackage{amscd}
\usepackage{verbatim}
\usepackage{mathrsfs}
\usepackage{color}

\allowdisplaybreaks[4]

\newcommand{\mf}{\mathfrak}
\newcommand{\g}{\mf{g}}

\newcommand{\gl}{\mf{gl}}

\allowdisplaybreaks[4]

\newcommand{\Z}{{\mathbb Z}}
\newcommand{\C}{{\mathbb C}}
\newcommand{\N}{{\mathbb N}}

\newcommand{\supp}{{\operatorname{Supp}}\xspace}

\renewcommand{\phi}{\varphi}

\renewcommand{\geq}{\geqslant}

\def\sl{\mathfrak{sl}}
\def\gl{\mathfrak{gl}}

\def\sl{\mathfrak{sl}}

\newtheorem{theorem}{Theorem}[section]

\newtheorem{lemma}[theorem]{Lemma}

\theoremstyle{remark}

\numberwithin{equation}{section}

\def\span{\mathrm{span}}
\def\supp{\mathrm{Supp}}

\newcommand\pxi[1]{\frac{\partial}{\partial\xi_{#1}}}
\newcommand\pt[1]{\frac{\partial}{\partial t_{#1}}}

\newenvironment{claim}[1]{\par\noindent\underline{Claim.}\space#1}{}

\begin{document}




\title[]{Tensor modules over Witt superalgebras}
\author{Yaohui Xue, Yan Wang}
\maketitle

\begin{abstract} In this paper, we study the tensor module $P\otimes M$ over the Witt superalgebra $W_{m,n}^+$ (resp. $W_{m,n}$), where $P$ is a simple module
over the Weyl superalgebra $K_{m,n}^+$ (resp. $K_{m,n}$) and $M$ is simple weight module over the general linear Lie superalgebra $\gl(m,n)$. We obtain the necessary and sufficient conditions for $P\otimes M$ to be simple, and determine all simple subquotient of $P\otimes M$ when it is not simple. All the work leads to completion of some classification problems on the weight representation theory of $W_{m,n}^+$ and $W_{m,n}$.

\vspace{0.3cm}
\noindent{\it Keywords}: Witt superalgebra, simple module, tensor module

\noindent{\it  MSC2020}: 17B10, 17B65, 17B66
\end{abstract}
%
%



\section{Introduction}
The notations $\Z, \Z_+, \N$ and $\C$ denote the sets of all integers, non-negative integers, positive integers and complex numbers, respectively.
Let $m,n\in \Z_+$ and at least one of them be nonzero. Let $A_{m,n}^+$ (resp. $A_{m,n}$) be the tensor superalgebra of the polynomial algebra $\C[t_1,\dots,t_m]$ (resp. Laurent polynomial algebra $\C[t_1^{\pm 1},\dots,t_m^{\pm 1}]$) in $m$ even variables $t_1,\dots,t_m$ and the exterior algebra in $n$ odd variables $\xi_1,\dots,\xi_n$. Omit $\otimes$ in $A_{m,n}^+$ and $A_{m,n}$ for convenience. Denote by $W_{m,n}^+$ (resp. $W_{m,n}$) the Lie superalgebra of super-derivations of $A_{m,n}^+$ (resp. $A_{m,n}$), which is called the Witt superalgebra.

There are many results in the weight representation theory of $W_{m,n}^+$ and $W_{m,n}$. The Harish-Chandra modules (weight modules with finite-dimensional weight spaces) over $W_{1,0}^+$ were classified by O. Mathieu in \cite{Ma}. A complete description of the supports of all simple weight modules over $W_{m,0}^+$ was given in \cite{PS}. The simple bounded modules (weight modules with bounded weight multiplicities) over $W_{m,0}^+$ are classified in \cite{XL1}. It is well known that the simple Harish-Chandra modules for the Virasoro algebra (the universal central extension of $W_{1,0}$) were conjectured by V. Kac and classified by O. Mathieu in \cite{Ma}, see \cite{Su2} for another approach. Y. Billig and V. Futorny completed the classification of the simple Harish-Chandra modules for $W_{m,0}$ in \cite{BF1}. Simple Harish-Chandra modules over $W_{0,n}$ were classified in \cite{DMP}. The simple strong Harish-Chandra modules (weight modules that are Harish-Chandra modules with respect to the Cartan subalgebra of $W_{m,0}$) over the $N=2$ Ramond algebra (the central extension of $W_{1,1}$) were classified in \cite{Liu1}. For more related results, please refer to \cite{BF2, CG, E1, E2, LZ2, Sh, Su1}.

For us, the most important studies in the weight representation theory of $W_{m,n}^+$ and $W_{m,n}$ are the following. In \cite{LX}, the simple bounded modules over $W_{m,n}^+$ were classified. Every such module is a simple quotient of a tensor module. In \cite{XL2}, the simple strong Harish-Chandra modules over $W_{m,n}$ were classified, which was also solved in \cite{BFIK}. Every such module is either a simple quotient of a tensor module or a module of highest weight type. For consistency, we denote the tensor modules that appear in the above two papers as $F(P, M)$ and let $F(P, M)=P\otimes M$, where $P$ is a simple module over the Weyl superalgebra $K_{m,n}^+$ (resp. $K_{m,n}$) and $M$ is a simple module over the general linear Lie superalgebra $\gl(m,n)$. It can be seen that the study of the tensor module $F(P,M)$ is very important for making the weight representation theory of $W_{m,n}^+$ and $W_{m,n}$ complete. In \cite{LLZ}, the tensor module over $W_{m,0}^+$ (resp. $W_{m,0}$) has been studied. The main difference, as well as the main difficult, is that there exist infinite cases for non-simple tensor module when $n>0$, while there exist only finite cases for non-simple tensor module when $n=0$. The special case, that $m=0$ and $M$ is a simple highest weight $\gl_n$-module, has been studied in \cite{S} via a different approach.

This paper is arranged as follows. In Section 2, we give some basic notations and results for our study. In Section 3, we prove the necessary conditions for the $W_{m,n}^+$-modules (resp. $W_{m,n}$-module) $F(P, M)$ to be not simple in Theorem 3.5. In Section 4, we determine all the simple subquotient in Theorem 4.8 under the necessary conditions in Theorem 3.5.

\section{Preliminaries}
All vector spaces and algebras in this paper are over $\C$. A super vector space $V$ is a vector space endowed with a $\mathbb{Z}_2$-gradation
$V=V_{\bar{0}}\oplus V_{\bar{1}}$. The parity of a homogeneous element $v\in V_{\bar{i}}$ is denoted by $|v|=\bar{i}\in \Z_2$. Throughout this paper, when we write $|v|$ for an element $v\in V$, we will always assume that $v$ is a homogeneous element.

Any module over a Lie superalgebra or an associative superalgebra is assumed to be $\Z_2$-graded.
A module $M$ over a Lie superalgebra or an associative superalgebra $\g$ is called simple if $M$ does not have proper submodules, and is called strictly simple if $M$ does not have $\g$-invariant subspaces except $0$ and $M$. Clearly, a strictly simple module must be simple. Define the parity-change $\Pi(M)$ of $M$ as follows: $\Pi(M)_{\bar 0}=M_{\bar 1}, \Pi(M)_{\bar 1}=M_{\bar 0}$, and the action of $\g$ on $\Pi(M)$ is the same as that of $\g$ on $M$.

\vspace{3mm}\noindent 2.1 {\bf Witt superalgebra.}
For any $\alpha=(\alpha_1,\dots,\alpha_m)\in\Z^m$ and $i_1,\dots,i_k\in\{1,\dots,n\}$, let $t^\alpha=t_1^{\alpha_1}\cdots t_m^{\alpha_m}$ and $\xi_{i_1,\dots,i_k}=\xi_{i_1}\cdots \xi_{i_k}$. Also, for any nonempty subset $I \subset \{1,\dots,n\}$, let $\xi_I=\xi_{l_1,\dots,l_k}$ with $l_1<\dots<l_k$ and $I=\{l_1,\dots,l_k\}$. In addition, set $\xi_\varnothing=1$.

Witt superalgebra $W_{m,n}^+$ (resp. $W_{m,n}$) has a standard basis
$$\{t^\alpha\xi_I\pt{i},t^\alpha\xi_I\frac{\partial}{\partial\xi_j}\mid \alpha\in\Z_+^m(\mbox{resp.}\ \alpha\in\Z^m),I\subset\{1,\dots,n\},i\in\{1,\dots,m\},j\in\{1,\dots,n\}\}$$
with the bracket defined by
$$[t^\alpha\xi_{I}\partial,t^{\alpha'}\xi_{I'}\partial']=t^\alpha\xi_I\partial(t^{\alpha'}\xi_{I'})\partial'
-(-1)^{(|I|+|\partial|)(|I'|+|\partial'|)}t^{\alpha'}\xi_{I'}\partial'(t^\alpha\xi_I)\partial,$$
where $\partial,\partial'\in\{\pt{1},\dots,\pt{m},\pxi 1,\dots,\pxi n\}.$

\vspace{3mm}\noindent 2.2 {\bf The general linear Lie superalgebra.}
Let $\gl(m,n)=\gl(m,n)_{\bar 0}\oplus \gl(m,n)_{\bar 1}$ be the general linear Lie superalgebra consisting of all $(m+n)\times(m+n)$ matrices with
\begin{eqnarray*}
&\gl(m,n)_{\bar 0}=\span\{E_{i,j}\mid i,j\in\{1,\dots,m\}\ \mbox{or} \ i,j\in\{m+1,\dots,m+n\}\},\\
&\gl(m,n)_{\bar 1}=\span\{E_{i,m+j},E_{m+j,i}\mid i\in\{1,\dots,m\},j\in\{1,\dots,n\}\},
\end{eqnarray*}
where $E_{i,j},i,j\in\{1,2,\dots,m+n\}$ is the $(i,j)$-th matrix unit. Also,
$\gl(m,n)$ has a $\Z$-gradation, i.e. $\gl(m,n)=\gl(m,n)_{-1}\oplus\gl(m,n)_0\oplus\gl(m,n)_1$, where
\begin{eqnarray*}
&\gl(m,n)_{-1}=\span\{E_{m+j,i}\mid i\in\{1,\dots,m\},j\in\{1,\dots,n\}\},\\
&\gl(m,n)_1=\span\{E_{i,m+j}\mid i\in\{1,\dots,m\},j\in\{1,\dots,n\}\}
\end{eqnarray*}
and $\gl(m,n)_0=\gl(m,n)_{\bar 0}$. Obviously, this $\Z$-gradation is consistent with the $\Z_2$-gradation of $\gl(m,n)$.

A $\gl(m,n)$-module $M$ is called a weight module if $M=\oplus_{\lambda\in\C^m,\mu\in\C^n}M_{(\lambda,\mu)}$, where
$$M_{(\lambda,\mu)}=\{v\in M\mid E_{i,i}(v)=\lambda_iv,E_{m+j,m+j}(v)=\mu_jv,i\in\{1,\dots,m\},j\in\{1,\dots,n\}\}$$
is called the weight space with weight $(\lambda,\mu)$.
Denote by
$$\supp(M)=\{(\lambda,\mu)\in\C^{m+n}\mid M_{(\lambda,\mu)}\neq 0\}$$
the support set of $M$.

Let $V$ be a module over the Lie algebra $\gl(m,n)_0$. Then $V$ could be viewed as a module over the Lie superalgebra $\gl(m,n)_0$ with $V_{\bar 0}=V$.
By extending $V$ trivially to a $\gl(m,n)_0\oplus\gl(m,n)_1$-module, we get the induced module $\mbox{Ind}_{\gl(m,n)_0\oplus\gl(m,n)_1}^{\gl(m,n)}(V)$,
which is called the Kac module $K(V)$ of $V$. It is easy to see that $K(V)$ is isomorphic to $\Lambda(\gl(m,n)_{-1})\otimes V$ as superspaces.

\begin{lemma}\cite[Theorem 4.1]{CM}\label{L(V)}
For any simple $\gl(m,n)_0$-module $V$, the module $K(V)$ has a unique maximal submodule and the unique simple top of $K(V)$ is denoted by $L(V)$. Any simple $\gl(m,n)$-module is isomorphic to $L(V)$ for some simple $\gl(m,n)_0$-module $V$ up to a parity-change.
\end{lemma}

Let $e_1,\dots,e_m$ be the standard basis of $\Z^m$ and $\C^m$. Then the general linear Lie algebra $\gl_m$ has a natural representation on $\C^m$ defined by $E_{i,j}e_k=\delta_{j,k}e_i$. For any $r\in\N$, let $V(r)$ be the $\gl_m$-module that is isomorphic to the exterior product $\Lambda^r(\C^m)=\C^m\wedge\dots\wedge\C^m$ ($r$ times) as vector spaces with the action given by
$$x(v_1\wedge\dots\wedge v_r)=\sum_{i=1}^rv_1\wedge\dots\wedge v_{i-1}\wedge xv_i\wedge\dots\wedge v_r,\forall x\in\gl_m.$$
Let $V(0)=\C$ be the trivial $\gl_m$-module and $v\wedge c=cv$ for any $v\in\C^m,c\in\C$.

\vspace{3mm}\noindent 2.3 {\bf Weyl superalgebra.}
The Weyl superalgebra $K_{m,n}^+$ is the simple associative superalgebra
$$\C[t_1,\dots,t_m,\xi_1,\dots,\xi_n,\pt{1},\dots,\pt{m},\pxi 1,\dots,\pxi n],$$
while $K_{m,n}$ is the simple associative superalgebra
$$\C[t_1^{\pm 1},\dots,t_m^{\pm 1},\xi_1,\dots,\xi_n,\pt{1},\dots,\pt{m},\pxi 1,\dots,\pxi n].$$

A weight module $V$ over $K_{m,n}^+$ (resp. $K_{m,n}$) is a module such that $t_1\pt{1},\dots,t_m\pt{m}$ and $\xi_1\pxi{1},\dots,\xi_n\pxi{n}$ act on $V$ diagonally.
Note that $K_{m,n}^+$ (resp. $K_{m,n}$) is a tensor product of $K_{(1,0)}^+,\dots,K_{(m,0)}^+,K_{(0,1)}^+\dots,K_{(0,n)}^+$ (resp. $K_{(1,0)},\dots,K_{(m,0)},K_{(0,1)}\dots,K_{(0,n)}$), where $K_{(i,0)}^+$ (resp. $K_{(i,0)}$) is the subalgebra of $K_{m,n}^+$ (resp. $K_{m,n}$) generated by $t_i,\frac{\partial}{\partial t_i}$ (resp. $t_i^{\pm 1},\pt{i}$) for $i\in\{1,\dots,m\}$, and $K_{(0,j)}^+$ (resp. $K_{(0,j)}$) is the subalgebra of $K_{m,n}^+$ (resp. $K_{m,n}$) generated by $\xi_j,\frac{\partial}{\partial \xi_j}$ for $j\in\{1,\dots,n\}$. For $i\in\{1,\dots,m\}$, $K_{(i,0)}^+\cong K_{1,0}^+$ and any simple weight module over
$K_{1,0}^+$ is one of the following (see for example Corollary 2.9 in \cite{GS})
\begin{equation}\label{+1}
t_1^{\lambda_1}\C[t_1^{\pm 1}], \lambda_1\notin\Z;\ \C[t_1];\ \C[t_1^{\pm 1}]/\C[t_1].
\end{equation}
For $i\in\{1,\dots,m\}$, $K_{(i,0)}\cong K_{1,0}$ and any simple weight module over $K_{1,0}$ is isomorphic to $t_1^{\lambda_1}\C[t_1^{\pm 1}]$ for some $\lambda_1\in\C$.
For $j\in\{1,\dots,n\}$, $K_{(0,j)}^+=K_{(0,j)}\cong K_{0,1}^+=K_{0,1}$ and any simple module over $K_{0,1}$ is isomorphic to $\C[\xi_1]$ up to a parity-change.
Let $V$ be a simple weight module over $K_{m,n}^+$ (resp. $K_{m,n}$). It is easy to see
$$V\cong V_{(1,0)}\otimes\dots \otimes V_{(m,0)}\otimes V_{(0,1)}\otimes\dots\otimes V_{(0,n)},$$
where each $V_{(i,0)}$ is a simple weight module over $K_{1,0}^+$ (resp. $K_{1,0}$) for $i\in \{1,\dots,m\}$ and each $V_{(0,j)}$ is a simple weight module over $K_{0,1}$
for $j\in \{1,\dots,n\}$. For some properties on simple weight $K_{m,n}^+$-modules (resp. $K_{m,n}$-modules), please refer to \cite{LX} and \cite{XL2}.

\vspace{3mm}\noindent 2.4 {\bf Simple weight $\gl_n$-module.}
For a given $\lambda\in\C^n$, any simple weight $K_{n,0}^+$-module $V$ with $\lambda\in\supp(V)$ is uniquely determined, denoted by $W(\lambda)$. Obviously, $\supp(W(\lambda))=X_1\times X_2\times\dots\times X_n\subset \lambda+\Z^n$ where $X_j=\lambda_j+\Z$ if $\lambda_j\notin\Z$, $X_j=\Z_+$ if $\lambda_j\in\Z_+$ and $X_j=- \N$ if $\lambda_j\in - \N$ by (\ref{+1}). Then $W(\lambda)$ has a basis $\{y(\lambda')\ |\ \lambda'\in \supp(W(\lambda))\}$ with
\begin{equation}\label{Wlam}
t_j\cdot y(\lambda')=y(\lambda'+e_j),\ \pt{j}\cdot y(\lambda')=\lambda'_jy(\lambda'-e_j),\forall j\in\{1,\dots,n\},
\end{equation}
where $y(\lambda')=0$ if $\lambda'\notin\supp(W(\lambda))$.

It is known that the Lie algebra $\gl_n$ can be identified with a Lie subalgebra of $K_{n,0}^+$ with $E_{l,j}=t_l\pt{j}$ for $l, j\in\{1,\dots,n\}$.
Let $\lambda\in\C^n$, define $|\lambda|=\sum_{j=1}^n\lambda_j$ and
$$N(\lambda)=\span\{y(\lambda')\mid \lambda'\in\supp(W(\lambda)),|\lambda'|=|\lambda|\}.$$
Notations $W(\lambda)$ and $N(\lambda)$, and the following two lemmas are both from \cite{BBL}.
\begin{lemma}\cite[Proposition 2.12]{BBL}
For $\lambda\in\C^n$, the $N(\lambda)$ is a simple $\gl_n$-module and a simple $\sl_n$-module.
\end{lemma}
\begin{lemma}\cite[Theorem 5.8]{BBL}\label{weight1dim}
Let $V$ be an infinite-dimensional simple weight $\sl_n$-module having all weight spaces one-dimensional. Then $V$ is isomorphic to $N(\lambda)$ for some $\lambda\in\C^n$.
\end{lemma}

The following lemma plays a key role in our method.
\begin{lemma}\label{gln}
Suppose that $n\geqslant 2$.
Let $V$ be a simple weight $\gl_n$-module satisfying that
$$(E_{l,l}-E_{l,j}E_{j,l}+E_{j,j}E_{l,l}) V=0,\forall l,j\in\{1,\dots,n\},l\neq j.$$
Then $V\cong N(\lambda)$ for some $\lambda\in\C^n$.
\end{lemma}
\begin{proof}
Suppose that $V$ is finite-dimensional and that $v\in V$ is a highest weight vector of weight $\mu$. For any $l,j\in\{1,\dots,n\}$ with $l>j$,
$$0=(E_{l,l}-E_{l,j}E_{j,l}+E_{j,j}E_{l,l}) v=\mu_l(1+\mu_j)v.$$
Since $V$ is finite-dimensional, we have $\mu_1-\mu_2,\dots,\mu_{n-1}-\mu_n\in\Z_+$. Then it is easy to know that $\mu=(p,0,\dots,0)$ or $(-1,\dots,-1,-p-1)$ for some $p\in\Z_+$. Note that the highest weight of $N(\mu)$ is $\mu$ by (\ref{Wlam}). Therefore, $V$ is just the $\gl_n$-module $N(\mu)$.

Now suppose that $V$ is an infinite-dimensional simple weight $\gl_n$-module. Let $\lambda\in\supp(V)$ and $v\in V_\lambda\setminus\{0\}$, then $U(\gl_n)v=V$ by the simplicity of $V$. Note that $U(\gl_n)$ is itself a weight $\gl_n$-module. So $U(\gl_n)_0v=V_\lambda$ and $U(\gl_n)_0$ is spanned by elements of the form $E_{l_1,j_1}\cdots E_{l_k,j_k}$ satisfying
$$e_{l_1}+\dots+e_{l_k}=e_{j_1}+\dots+e_{j_k}.$$

\begin{claim}
$E_{l_1,j_1}\cdots E_{l_k,j_k} v\in\C v$ for all $E_{l_1,j_1}\cdots E_{l_k,j_k}\in U(\gl_n)_0$.
\end{claim}

We will prove this by induction on $k$. This is true for $k=1$ or $2$ since $v$ is a weight vector and
$$E_{l,j}E_{j,l} v=(E_{l,l}+E_{j,j}E_{l,l}) v,\forall l\neq j.$$
For $k\geqslant 3$, suppse $l_r\neq j_r$ for $r\in\{1,\dots,k\}$. Let $s\geq 2$ be the minimal integer such that $l_1=j_s$ or $l_s=j_1$, then
$$E_{l_1,j_1}\cdots E_{l_k,j_k} v=E_{l_2,j_2}\cdots E_{l_{s-1},j_{s-1}}E_{l_1,j_1}E_{l_s,j_s}\cdots E_{l_k,j_k} v.$$

If $l_1=j_s$ and $l_s=j_1$, we have
$$E_{l_2,j_2}\cdots E_{l_{s-1},j_{s-1}}E_{l_1,j_1}E_{l_s,j_s}\cdots E_{l_k,j_k} v\in E_{l_2,j_2}\cdots E_{l_{s-1},j_{s-1}}E_{l_{s+1},j_{s+1}}\cdots E_{l_k,j_k}(\C v),$$
which is contained in $\C v$ by the induction hypothesis.

If $l_1\neq j_s$ and $l_s=j_1$ (the case that $l_1=j_s$ and $l_s\neq j_1$ could be done similarly), we have
\begin{equation*}\begin{split}
E_{l_1,j_1}E_{l_s,j_s}&=E_{l_1,j_1}[E_{j_1,l_1},E_{l_1,j_s}]\\
&=E_{l_1,j_1}E_{j_1,l_1}E_{l_1,j_s}-E_{l_1,j_1}E_{l_1,j_s}E_{j_1,l_1}\\
&=E_{l_1,j_1}E_{j_1,l_1}E_{l_1,j_s}-E_{l_1,j_s}E_{l_1,j_1}E_{j_1,l_1}.
\end{split}\end{equation*}
Then
$$E_{l_2,j_2}\cdots E_{l_{s-1},j_{s-1}}E_{l_1,j_1}E_{l_s,j_s}\cdots E_{l_k,j_k} v\in E_{l_2,j_2}\cdots E_{l_{s-1},j_{s-1}}E_{l_1,j_s}E_{l_{s+1},j_{s+1}}\cdots E_{l_k,j_k}(\C v).$$
The Claim follows from the induction hypothesis.

The Claim implies that $V$ is a weight module having all weight spaces one-dimensional. Since $V$ is a simple $\gl_n$-module, $E_{1,1}+\dots+E_{n,n}$ acts on $V$ as a scalar. Then $V$ is also a weight $\sl_n$-module with one-dimensional weight spaces. By Lemma \ref{weight1dim}, the $\gl_n$-module $V=N(\lambda)\otimes \C v_b$, where $\lambda\in\C^n$ and $\C v_b$ is the one-dimensional $\gl_n$-module with $\sl_n v_b=0$ and $(E_{1,1}+\dots+E_{n,n}) v_b=bv_b$ for $b\in \C$.
Then we have
$$E_{l,l}\cdot (y(\lambda')\otimes v_b)=(\lambda'_l+\frac{b}{n})y(\lambda')\otimes v_b,l=1,\dots,n.$$
Furthermore, for any $l,j\in\{1,\dots,n\}$ with $l\neq j$ and $\lambda'\in\supp(W(\lambda))$,
\begin{equation*}\begin{split}
0=&(E_{l,l}-E_{l,j}E_{j,l}+E_{j,j}E_{l,l})\cdot (y(\lambda')\otimes v_b)\\
=&((\lambda'_l+\frac{b}{n})-\lambda'_l(\lambda'_j+1)+(\lambda'_l+\frac{b}{n})(\lambda'_j+\frac{b}{n})) (y(\lambda')\otimes v_b)\\
=&\frac{b}{n}(\lambda'_l+\lambda'_j+\frac{b}{n}+1)(y(\lambda')\otimes v_b).
\end{split}\end{equation*}
Then
\begin{equation}\label{b}
b=0\ \text{or}\ -n\lambda'_l-n\lambda'_j-n.
\end{equation}

If $b=0$, $N(\lambda)\otimes \C v_b$ is just $N(\lambda)$. Suppose $b\neq 0$, then $b=-n\lambda'_l-n\lambda'_j-n$.
If $n\geqslant 3$, there is a $\lambda''\in\supp(W(\lambda))$ such that $\lambda''_l=\lambda'_l,\lambda''_j\neq\lambda'_j$ or $\lambda''_j=\lambda'_j,\lambda''_l\neq\lambda'_l$. By (\ref{b}), there is $-n\lambda''_l-n\lambda''_j-n=b=-n\lambda'_l-n\lambda'_j-n$, which is a contradiction. So $n=2$ and $b=-2|\lambda|-2$. It is easy to verify that the map $\tau:N(\lambda)\otimes\C v_b\rightarrow N((-\lambda_2-1,-\lambda_1-1))$, defined by
\begin{eqnarray*}
&\tau(y((\lambda_1+s,\lambda_2-s))\otimes v_b)=(-1)^s\binom{\lambda_1+s}{s}/\binom{\lambda_2}{s}y((-\lambda_2+s-1,-\lambda_1-s-1)), \\
&\tau(y((\lambda_1-s,\lambda_2+s))\otimes v_b)=(-1)^s\binom{\lambda_2+s}{s}/\binom{\lambda_1}{s}y((-\lambda_2-s-1,-\lambda_1+s-1)),\\
&(\lambda_1+s, \lambda_2-s)\in\supp(W(\lambda)),(\lambda_1-s, \lambda_2+s)\in\supp(W(\lambda)),s\in\Z_+,
\end{eqnarray*}
is an isomorphism of $\gl_2$-modules. Hence $V\cong N((-\lambda_2-1,-\lambda_1-1))$.
\end{proof}

\section{Simplicity of $F(P,M)$}
From \cite{LX}, we know that there is a Lie superalgebra homomorphism $\pi$ from $W_{m,n}^+$ to the tensor superalgebra $K_{m,n}^+\otimes U(\gl(m,n))$ given by
\begin{equation}\label{pi}\begin{split}
\pi(t^\alpha\xi_I\pt{i})=&t^\alpha\xi_I\pt{i}\otimes 1+\sum_{s=1}^m\pt{s}(t^\alpha\xi_I)\otimes E_{s,i}+(-1)^{|I|-1}\sum_{l=1}^n\pxi l(t^\alpha\xi_I)\otimes E_{m+l,i},\\
\pi(t^\alpha\xi_I\pxi j)=&t^\alpha\xi_I\pxi j\otimes 1+\sum_{s=1}^m\pt{s}(t^\alpha\xi_I)\otimes E_{s,m+j}+(-1)^{|I|-1}\sum_{l=1}^n\pxi l(t^\alpha\xi_I)\otimes E_{m+l,m+j},
\end{split}\end{equation}
where $\alpha\in\Z^m_+$, $I\subset\{1,\dots,n\}$, $i\in\{1,\dots,m\}$ and $j\in\{1,\dots,n\}$. It could be verified that (\ref{pi}), with $\alpha\in\Z^m$, also gives a Lie superalgebra homomorphism $\pi$ from $W_{m,n}$ to the tensor superalgebra $K_{m,n}\otimes U(\gl(m,n))$.

Let $P$ be a module over $K_{m,n}^+$ (resp. $K_{m,n}$) and $M$ be a $\gl(m,n)$-module. Then we have the tensor module $F(P,M)=P\otimes M$ over $K_{m,n}^+\otimes U(\gl(m,n))$ (resp. $K_{m,n}\otimes U(\gl(m,n))$). It follows that $F(P,M)$ is a module over $W_{m,n}^+$ (resp. $W_{m,n}$) with the action given by
$$x\cdot (u\otimes v)=\pi(x)(u\otimes v)$$
for any $x\in W_{m,n}^+ (\mbox{resp}.\ W_{m,n})$.

Now, we give a lemma that will be used in the following proof.
\begin{lemma}\cite[Lemma 2.2]{XL2}\label{density}
Let $B,B'$ be two unital associative superalgebras such that $B'$ has a countable basis, $R=B\otimes B'$.
\begin{itemize}
  \item[(1)]Let $M$ be a $B$-module and $M'$ be a strictly simple $B'$-module. Then $M\otimes M'$ is a simple $R$-module if and only if $M$ is a simple $B$-module.
  \item[(2)]Suppose that $V$ is a simple $R$-module and $V$ contains a strictly simple $B'=\C\otimes B'$-submodule $M'$. Then $V\cong M\otimes M'$ for some simple $B$-module $M$.
\end{itemize}
\end{lemma}
\begin{lemma}\label{V_1}
Let $P$ be a simple module over $K_{m,n}^+$ (resp. $K_{m,n}$) and $M$ be a simple $\gl(m,n)$-module. Suppose the $W_{m,n}^+$-module (resp. $W_{m,n}$-module) $F(P,M)$ is not simple, then $M$ is of the form $L(V_1\otimes V_2)$ up to a parity-change, where $V_1$ is one of the finite-dimensional simple $\gl_m$-modules $V(r),r\in\{0,\dots,m\}$, and $V_2$ is a simple $\gl_n$-module.
\end{lemma}
\begin{proof}
We will prove this only for $W_{m,n}^+$ and this proof is also valid for $W_{m,n}$.
Let $F'$ be a nonzero proper submodule of $F(P,M)$. Let $\sum_{p=1}^qu_p\otimes v_p$ be a nonzero homogeneous vector in $F'$, where $u_p$'s are homogeneous vectors in $P$ and $v_p$'s are homogeneous vectors in $M$.

\underline{Claim} 1. For any $x\in K_{m,n}^+$ and $i, k, r\in\{1,\dots,m\}$, we have
$$\sum_{p=1}^qxu_p\otimes(\delta_{r,i}E_{r,k}-E_{r,i}E_{r,k})v_p\in F'.$$

First, for any $i\in\{1,\dots,m\}$ and $j\in\{1,\dots,n\}$, note that
\begin{equation}\label{eq1}
\pt{i}\cdot(\sum_{p=1}^qu_p\otimes v_p)=\sum_{p=1}^q\pt{i}u_p\otimes v_p,\ \pxi{j}\cdot(\sum_{p=1}^qu_p\otimes v_p)=\sum_{p=1}^q\pxi{j}u_p\otimes v_p.
\end{equation}
Fix $i, k$ and $r$. For $d\in\{0,1,2\}$ and $\beta\in\Z^m_+$ with $\beta_r\geqslant 2$, we have
\begin{eqnarray*}
&&t_r^d\pt{i}\cdot t^{\beta-de_r}\xi_I\pt{k}\cdot(\sum_{p=1}^q u_p\otimes v_p)\\
&=&\sum_{p=1}^qt_r^d\pt{i}\cdot\big(t^{\beta-de_r}\xi_I\pt{k}u_p\otimes v_p+\sum_{s=1}^m(\beta_s-d\delta_{r,s})t^{\beta-de_r-e_s}\xi_Iu_p\otimes E_{s,k}v_p\\
&&+(-1)^{|I|-1}\sum_{l=1}^n(-1)^{|u_p|}t^{\beta-de_r}\pxi{l}(\xi_I)u_p\otimes E_{m+l,k}v_p\big)\\
&=&\sum_{p=1}^q\big(t_r^d\pt{i}t^{\beta-de_r}\xi_I\pt{k}u_p\otimes v_p+dt_r^{d-1}t^{\beta-de_r}\xi_I\pt{k}u_p\otimes E_{r,i}v_p\big)\\
&&+\sum_{p=1}^q\sum_{s=1}^m(\beta_s-d\delta_{r,s})\big(t_r^d\pt{i}t^{\beta-de_r-e_s}\xi_Iu_p\otimes E_{s,k}v_p+dt_r^{d-1}t^{\beta-de_r-e_s}\xi_Iu_p\otimes E_{r,i}E_{s,k}v_p\big)\\
&&+\sum_{p=1}^q\sum_{l=1}^n(-1)^{|I|-1+|u_p|}\big(t_r^d\pt{i}t^{\beta-de_r}\pxi{l}(\xi_I)u_p\otimes E_{m+l,k}v_p\\
&&+dt_r^{d-1}t^{\beta-de_r}\pxi{l}(\xi_I)u_p\otimes E_{r,i}E_{m+l,k}v_p\big)\\
&=&\sum_{p=1}^q\big(t^\beta\pt{i}\xi_I\pt{k}u_p\otimes v_p+(\beta_i-d\delta_{i,r})t^{\beta-e_i}\xi_I\pt{k}u_p\otimes v_p
+dt^{\beta-e_r}\xi_I\pt{k}u_p\otimes E_{r,i}v_p\big)\\
&&+\sum_{p=1}^q\sum_{s=1}^m(\beta_s-d\delta_{r,s})\big(t^{\beta-e_s}\pt{i}\xi_Iu_p\otimes E_{s,k}v_p+(\beta_i-d\delta_{r,i}-\delta_{s,i})t^{\beta-e_s-e_i}\xi_Iu_p\otimes E_{s,k}v_p\\
&&+dt^{\beta-e_r-e_s}\xi_Iu_p\otimes E_{r,i}E_{s,k}v_p\big)+\sum_{p=1}^q\sum_{l=1}^n(-1)^{|I|-1+|u_p|}\big(t^\beta\pt{i}\pxi{l}(\xi_I)u_p\otimes E_{m+l,k}v_p\\
&&+(\beta_i-d\delta_{i,r})t^{\beta-e_i}\pxi{l}(\xi_I)u_p\otimes E_{m+l,k}v_p+dt^{\beta-e_r}\pxi{l}(\xi_I)u_p\otimes E_{r,i}E_{m+l,k}v_p\big)\\
&=&d^2\sum_{p=1}^qt^{\beta-2e_r}\xi_Iu_p\otimes(\delta_{r,i}E_{r,k}-E_{r,i}E_{r,k})v_p+dx_1+x_0,
\end{eqnarray*}
where $x_1$ and $x_0$ are determined by $\beta,i,k,r$ and are independent of $d$. Since $d$ is arbitrary in $\{0,1,2\}$, it is easy to know that
\begin{equation}\label{eq2}
\sum_{p=1}^qt^\alpha\xi_Iu_p\otimes(\delta_{r,i}E_{r,k}-E_{r,i}E_{r,k})v_p\in F',\forall \alpha\in\Z^m_+.
\end{equation}
By the structure of $K_{m,n}^+$, Claim 1 follows easily from (\ref{eq1}) and (\ref{eq2}).

\underline{Claim} 2. If $u_1,\dots,u_q$ are linearly independent, then we have $(\delta_{r,i}E_{r,k}-E_{r,i}E_{r,k})v_p=0$ for any $i, k, r\in\{1,\dots,m\}$ and $p\in\{1,\dots,q\}$.

Since $P$ is a simple $K_{m,n}^+$-module, $P$ is a strictly simple module by Lemma 3.3 in \cite{LX}. By the Jacobson density theorem in ring theory, for any $u\in P$, there is an $x\in K_{m,n}^+$ such that $xu_p=\delta_{1,p}u,p\in\{1, \dots,q\}$. By Claim 1, we have
$$P\otimes(\delta_{r,i}E_{r,k}-E_{r,i}E_{r,k})v_p\subset F'.$$
Let $M'=\{v\in M\ |\ P\otimes v\subset F'\}$, then $M'$ is clearly $\Z_2$-graded. For any $i,s\in\{1,\dots,m\}, j,l\in\{1,\dots,n\}$ and $u\otimes v\in F'$, from
\begin{eqnarray*}
&t_i\pt{s}\cdot(u\otimes v)=t_i\pt{s}u\otimes v+u\otimes E_{i,s}v,\\
&\xi_j\pt{s}\cdot(u\otimes v)=\xi_j\pt{s}u\otimes v+(-1)^{|u|}u\otimes E_{m+j,s}v,\\
&t_i\pxi{l}\cdot(u\otimes v)=t_i\pxi{l}u\otimes v+(-1)^{|u|}u\otimes E_{i,m+l}v,\\
&\xi_j\pxi{l}\cdot(u\otimes v)=\xi_j\pxi{l}u\otimes v+u\otimes E_{m+j,m+l}v,
\end{eqnarray*}
we see that $M'$ is a $\gl(m,n)$-submodule of $M$. Moreover, $P\otimes M'\subset F'\subsetneq F(P,M)$ implies that $M'\neq M$. So, from the simplicity of $M$, we
get $M'=0$. Thus Claim 2 is true.

From now on we assume that $u_1,\dots,u_q$ are linearly independent.

\underline{Claim} 3. For any $i, k, r\in\{1,\dots,m\}$, we have $(\delta_{r,i}E_{r,k}-E_{r,i}E_{r,k}) M=0$.

For $i_1,i_2\in\{1,\dots,m\}$, we have
$$t_{i_1}\pt{i_2}\cdot(\sum_{p=1}^qu_p\otimes v_p)=\sum_{p=1}^qt_{i_1}\pt{i_2}u_p\otimes v_p+\sum_{p=1}^qu_p\otimes E_{i_1,i_2}v_p\in F'.$$
By Claims 1 and 2, for any $x\in K_{m,n}^+$, we have
\begin{eqnarray*}
&&\sum_{p=1}^qxt_{i_1}\pt{i_2}u_p\otimes(\delta_{r,i}E_{r,k}-E_{r,i}E_{r,k})v_p
+\sum_{p=1}^qxu_p\otimes(\delta_{r,i}E_{r,k}-E_{r,i}E_{r,k})E_{i_1,i_2}v_p\\
&=&\sum_{p=1}^qxu_p\otimes(\delta_{r,i}E_{r,k}-E_{r,i}E_{r,k})E_{i_1,i_2}v_p\in F'.
\end{eqnarray*}
Since $u_p$'s are linearly independent, we deduce that
$$P\otimes(\delta_{r,i}E_{r,k}-E_{r,i}E_{r,k})E_{i_1,i_2}v_p\subset F'.$$
Similarly, we have
$$P\otimes(\delta_{r,i}E_{r,k}-E_{r,i}E_{r,k})E_{l_1,l_2}v_p\subset F',\forall l_1,l_2\in\{1,\dots,m+n\}.$$
This means that
$$(\delta_{r,i}E_{r,k}-E_{r,i}E_{r,k})E_{l_1,l_2}v_p=0,\forall l_1,l_2\in\{1,\dots,m+n\}.$$
By repeatedly doing this procedure we deduce that
\begin{eqnarray*}
(\delta_{r,i}E_{r,k}-E_{r,i}E_{r,k})U(\gl(m,n))v_p=(\delta_{r,i}E_{r,k}-E_{r,i}E_{r,k})M=0.
\end{eqnarray*}
Thus Claim 3 is true.

By Lemma \ref{L(V)}, there is a simple $\gl(m,n)_0$-submodule $V$ of $M$ such that $M=L(V)$ up to a parity-change. Let $v$ be a nonzero vector in $V$. By the PBW Theorem, Claim 3 implies that the $\gl_m$-submodule $U(\gl_m)v$ of $V$ is finite-dimensional. Let $V_1$ be a simple $\gl_m$-submodule of $U(\gl_m)v$. Since
$$(E_{i,i}-1)E_{i,i} V_1=0,\forall i\in\{1,\dots,m\},$$
it is well-known that $V_1$ is isomorphic to one of the simple $\gl_m$-modules $V(r),r\in\{0,\dots,m\}$.

Note that $U(\gl(m,n)_0)\cong U(\gl_m)\otimes U(\gl_n)$. Hence, by Lemma \ref{density}, there is a simple $\gl_n$-module $V_2$ such that $V\cong V_1\otimes V_2$. Then $M\cong L(V_1\otimes V_2)$ or $\Pi(L(V_1\otimes V_2))$.
\end{proof}

Let $P$ be a simple module over $K_{m,n}^+$ (resp. $K_{m,n}$). It is known that any simple $K_{0,n}^+$-module is isomorphic to $\Lambda(n)$ up to a parity-change. Since $K_{m,n}^+\cong K_{m,0}^+\otimes K_{0,n}^+$ (resp. $K_{m,n}\cong K_{m,0}\otimes K_{0,n}$) and $K_{0,n}^+$ ($=K_{0,n}$) is finite-dimensional, there is
\begin{equation}\label{P}
P\cong P_1\otimes\Lambda(n)
\end{equation}
by Lemma \ref{density}, where $P_1$ is a simple module over $K_{m,0}^+$ (resp. $K_{m,0}$).
\begin{lemma}\label{P1oV}
Let $P$ be a simple module over $K_{m,n}^+$ (resp. $K_{m,n}$), $V_1$ be a simple $\gl_m$-module and $V_2$ be a non-trivial simple $\gl_n$-module. Then any nonzero submodule of $F(P,L(V_1\otimes V_2))$ contains $P_1\otimes(V_1\otimes V_2)$, where $P_1=\{u\in P\mid \pxi{j}u=0,j\in\{1,\dots,n\}\}$. Thus, $F(P,L(V_1\otimes V_2))$ has a unique simple submodule.
\end{lemma}
\begin{proof}
We will prove this only for $W_{m,n}^+$ and this proof is also valid for $W_{m,n}$. Let $V$ be the simple $\gl(m,n)_0$-module $V_1\otimes V_2$ and $F'$ be a nonzero submodule of $F(P,L(V))$. Let $\sum_{p=1}^qu_p\otimes v_p$ be a nonzero homogeneous vector in $F'$, where $u_p$'s are homogeneous vectors in $P$ and linearly independent, and $v_p$'s are homogeneous vectors in $L(V)$.

For any $j\in\{1,\dots,n\}$, there is $\pxi{j}\cdot(\sum_{p=1}^qu_p\otimes v_p)=\sum_{p=1}^q\pxi{j}u_p\otimes v_p$.
By (\ref{P}), without loss of generality, we could assume that $\sum_{p=1}^qu_p\otimes v_p\in F'\setminus\{0\}$ with
$$\pxi{j}u_p=0,\forall j\in\{1,\dots,n\},\ p\in\{1,\dots,q\}.$$
Namely, we can say that each $u_p\in P_1$.

By the proof of Theorem 4.1 in \cite{CM}, the $\gl(m,n)$-module $L(V)$ is $\Z$-graded with respect to the $\Z$-gradation of $\gl(m,n)$ with the top non-zero graded component being $V$. From the PBW Theorem and the simplicity of $L(V)$, for any non-zero $v\in L(V)$, there is an $x\in U(\gl(m,n)_1)$ such that $xv\in V\setminus\{0\}$. For any $s\in\{1,\dots,m\}$ and $j\in\{1,\dots,n\}$, we have
$$t_s\pxi{j}\cdot(\sum_{p=1}^qu_p\otimes v_p)=\sum_{p=1}^q(-1)^{|u_p|}u_p\otimes E_{s,m+j}v_p\in F'.$$
So
$$\sum_{p=1}^q(-1)^{|u_p||x|}u_p\otimes xv_p\in F',\forall x\in U(\gl(m,n)_1).$$
Thus, without loss of generality, we could also assume that $\sum_{p=1}^qu_p\otimes v_p\in F'\setminus\{0\}$ with $v_p\in V$. Then $\gl(m,n)_1v_p=0,p\in\{1,\dots,q\}$.

Since $V_2$ is a non-trivial simple $\gl_n$-module, there exist $l,j\in\{1,\dots,n\}$ such that $E_{m+l,m+j}v_1\neq 0$. For any $\alpha\in\Z_+^m$, we have
$$t^\alpha\xi_l\pxi{j}\cdot(\sum_{p=1}^qu_p\otimes v_p)=\sum_{p=1}^qt^\alpha u_p\otimes E_{m+l,m+j}v_p\in F'.$$
Note that $\pt{i}\cdot (u\otimes v)=\pt{i}u\otimes v$ for any $i\in\{1,\dots,m\},u\in P$ and $v\in L(V)$. It follows that $\sum_{p=1}^qxu_p\otimes E_{m+l,m+j}v_p\in F'$ for any $x\in K_{m,0}^+$. From (\ref{P}), $P_1$ is a simple $K_{m,0}^+$-submodule of $P$. So, by the Jacobson density theorem, for any $u\in P_1$, there is an $x\in K_{m,0}^+$ such that $xu_p=\delta_{1,p}u, p\in\{1,\dots, q\}$. Therefore $P_1\otimes E_{m+l,m+j}v_1\subset F'$.

For any $u\in P_1, s, i\in\{1,\dots,m\}$, we have
$$t_s\pt{i}\cdot(u\otimes E_{m+l,m+j}v_1)=t_s\pt{i}u\otimes E_{m+l,m+j}v_1+u\otimes E_{s,i}E_{m+l,m+j}v_1\in F'.$$
So $u\otimes E_{s,i}E_{m+l,m+j}v_1\in F'$. Also, for any $l',j'\in\{1,\dots,n\}$ we have
$$\xi_{l'}\pxi{j'}\cdot (u\otimes E_{m+l,m+j}v_1)=u\otimes E_{m+l',m+j'}E_{m+l,m+j}v_1\in F'.$$
Then $u\otimes U(\gl(m,n)_0)E_{m+l,m+j}v_1=u\otimes V\subset F'$. Hence, $P_1\otimes V\subset F'$.
\end{proof}

\begin{lemma}\label{0m-1}
Let $P$ be a simple module over $K_{m,n}^+$ (resp. $K_{m,n}$), $V_1$ be one of the simple $\gl_m$-modules $V(r),r\in\{0,\dots,m-1\}$ and $V_2$ be a non-trivial simple $\gl_n$-module. Then the $W_{m,n}^+$-module (resp. $W_{m,n}$-module) $F(P,L(V_1\otimes V_2))$ is simple.
\end{lemma}
\begin{proof}
We will prove this only for $W_{m,n}^+$ and this proof is also valid for $W_{m,n}$. Let $V$ be the simple $\gl(m,n)_0$-module $V_1\otimes V_2$. Suppose that $F(P,L(V))$ is not simple and $F'$ is a nonzero proper submodule of $F(P,L(V))$. By Lemma \ref{P1oV}, $P_1\otimes V\subset F'$, where
$P_1=\{u\in P\mid \pxi{j}u=0,j\in\{1,\dots,n\}\}$ is a simple $K_{m,0}^+$-submodule of $P$. Since $V_1$ is one of the simple $\gl_m$-modules $V(r),r\in\{0,\dots,m-1\}$, there is a nonzero vector $v$ in $V$ with $E_{i,i}v=0$ for some $i\in\{1,\dots,m\}$.

By Claim 3 in the proof of Lemma \ref{V_1}, $(E_{i,i}-1)E_{i,i}L(V)=0, \forall i\in\{1,\dots,m\}$, which implies that $E_{i,i}$ only has eigenvalues $0$ and $1$ in $L(V)$. So, for any $l\in\{1,\dots,n\}$, $E_{i,i}E_{m+l,i}v=-E_{m+l,i}v$ implies that $E_{m+l,i}v=0$. Let $I_n=\{1,\dots,n\}$, then
$$\xi_{I_n}\pt{i}\cdot (u\otimes v)=\xi_{I_n}\pt{i}(u)\otimes v\in F',\forall u\in P_1.$$
Also,
$$t_i\xi_{I_n}\pt{i}\cdot (u\otimes v)=t_i\xi_{I_n}\pt{i}(u)\otimes v\in F', \forall u\in P_1.$$
Therefore,
$$\xi_{I_n}(u)\otimes v=\xi_{I_n}(\pt{i}t_i-t_i\pt{i})(u)\otimes v\in F', \forall u\in P_1.$$
So from $\xi_{I_n}(P_1)\otimes v\subset F'$, we can go further and say that $P\otimes v\subset F'$.
Note that $\{w\in L(V)\mid P\otimes w\subset F'\}$ is a $\gl(m,n)$-submodule of $L(V)$. By the simplicity of $L(V)$, $\{w\in L(V)\mid P\otimes w\subset F'\}=L(V)$. Consequently, $F'=F(P,L(V))$, which is a contradiction. Hence, $F(P,L(V))$ is simple.
\end{proof}

\begin{lemma}\label{V_2}
Let $P$ be a simple module over $K_{m,n}^+$ (resp. $K_{m,n}$), $V_1$ be a simple $\gl_m$-module and $V_2$ be a simple weight $\gl_n$-module. Suppose that the $W_{m,n}^+$-module (resp. $W_{m,n}$-module) $F(P,L(V_1\otimes V_2))$ is not simple. Then $V_2$ is the simple $\gl_n$-module $N(\lambda)$ for some $\lambda\in\C^n$.
\end{lemma}
\begin{proof}
We will prove this only for $W_{m,n}^+$ and this proof is also valid for $W_{m,n}$. This lemma is trivial if $n=1$. Suppose that $n\geqslant 2$.

Let $V$ be the simple $\gl(m,n)_0$-module $V_1\otimes V_2$ and $F'$ be a nonzero proper submodule of $F(P,L(V))$. Let $\sum_{p=1}^qu_p\otimes v_p$ be a nonzero homogeneous vector in $F'$, where $u_p$'s are homogeneous vectors in $P$ and $v_p$'s are homogeneous vectors in $L(V)$. By the proof of Lemma \ref{P1oV}, we can assume that
$\pxi{j}u_p=0$ and $v_p\in V$ for $j\in\{1,\dots,n\}$ and $p\in\{1,\dots,q\}$.

For $l',j,j'\in\{1,\dots,n\}$ with $l'\neq j', \alpha\in\Z^m_+$ and $I_n=\{1,\dots,n\}$, we have
\begin{eqnarray*}
&&(-1)^{n-1}\xi_{l'}\xi_{j'}\pxi{j'}\cdot t^\alpha\xi_{I_n}\pxi{j}\cdot(\sum_{p=1}^qu_p\otimes v_p)\\
&=&\xi_{l'}\xi_{j'}\pxi{j'}\cdot\sum_{p=1}^q\sum_{l=1}^n\pxi{l}(t^\alpha\xi_{I_n})u_p\otimes E_{m+l,m+j}v_p\\
&=&\sum_{p=1}^q\sum_{l=1}^n\big (\xi_{l'}\xi_{j'}\pxi{j'}\pxi{l}(t^\alpha\xi_{I_n})u_p\otimes E_{m+l,m+j}v_p\\
&&-\xi_{j'}\pxi{l}(t^\alpha\xi_{I_n})u_p\otimes E_{m+l',m+j'}E_{m+l,m+j}v_p+\xi_{l'}\pxi{l}(t^\alpha\xi_{I_n})u_p\otimes E_{m+j',m+j'}E_{m+l,m+j}v_p\big )\\
&=&\sum_{p=1}^q \big (\xi_{l'}t^\alpha\pxi{l'}(\xi_{I_n})u_p\otimes E_{m+l',m+j}v_p
-\xi_{j'}\pxi{j'}(t^\alpha\xi_{I_n})u_p\otimes E_{m+l',m+j'}E_{m+j',m+j}v_p\\
&&+\xi_{l'}\pxi{l'}(t^\alpha\xi_{I_n})u_p\otimes E_{m+j',m+j'}E_{m+l',m+j}v_p\big )\\
&=&\sum_{p=1}^q\big (t^\alpha\xi_{I_n}u_p\otimes(E_{m+l',m+j}-E_{m+l',m+j'}E_{m+j',m+j}+E_{m+j',m+j'}E_{m+l',m+j})v_p\big )\in F'.
\end{eqnarray*}
This result, combined with (\ref{eq1}), leads to the following conclusion
$$\sum_{p=1}^q(xu_p\otimes(E_{m+l',m+j}-E_{m+l',m+j'}E_{m+j',m+j}+E_{m+j',m+j'}E_{m+l',m+j})v_p)\in F',\forall x\in K_{m,n}^+.$$
As we did in the proof of Lemma \ref{V_1}, we could prove that
$$(E_{m+l',m+j}-E_{m+l',m+j'}E_{m+j',m+j}+E_{m+j',m+j'}E_{m+l',m+j})V=0$$
for all $l',j,j'\in\{1,\dots,n\}$ with $l'\neq j'$. By letting $l'=j$, we have
$$(E_{m+j,m+j}-E_{m+j,m+j'}E_{m+j',m+j}+E_{m+j',m+j'}E_{m+j,m+j})V=0.$$
Therefore, the $\gl_n$-module $V_2$ satisfies the following condition
$$(E_{j,j}-E_{j,j'}E_{j',j}+E_{j',j'}E_{j,j})V_2=0,\forall j,j'\in\{1,\dots,n\},j\neq j'.$$
By Lemma \ref{gln}, $V_2=N(\lambda)$ for some $\lambda\in\C^n$.
\end{proof}

\begin{theorem}\label{notsimple}
Let $P$ be a simple module over $K_{m,n}^+$ (resp. $K_{m,n}$) and $M$ be a simple weight $\gl(m,n)$-module. Suppose the $W_{m,n}^+$-module (resp. $W_{m,n}$-module) $F(P,M)$ is not simple, then $M$ is of the form $L(V(r)\otimes N(\lambda))$ for $r\in\{0,\dots, m\}$ and $\lambda\in\C^n$ up to a parity-change, where $\lambda={\bf 0}\in \C^n$ if $r\neq m$.
\end{theorem}
\begin{proof}
This theorem follows from Lemmas \ref{V_1}, \ref{0m-1} and \ref{V_2}.
\end{proof}

\section{Subquotients of $F(P,M)$}

For $\lambda\in\C^n$, denote by $M(\lambda)$ the $\gl(m,n)_0$-module that is a tensor product of the $\gl_m$-module $\sum_{r=0}^mV(r)$ and the $\gl_n$-module $W(\lambda)$. Omit the $\otimes$ in $M(\lambda)$ for convenience. For $i\in\{0,1\}$, let
\begin{eqnarray*}
&M(\lambda)_{\bar i}=\span\{e_{i_1}\wedge\dots\wedge e_{i_r}y(\lambda')\mid r\in\{0,\dots,m\},i_1,\dots,i_r\in\{1,\dots,m\},\\
&\lambda'\in\supp(W(\lambda)),\overline{|\lambda-\lambda'|}=\bar i\in\Z_2\}.
\end{eqnarray*}

The proof of the following lemma is given in detail in Appendix.
\begin{lemma}\label{4.1}
$M(\lambda)$ is a $\gl(m,n)$-module with the actions of $\gl(m,n)_{\pm 1}$ defined as follows:
\begin{equation*}\begin{split}
&E_{i,m+j}\cdot (e_{i_1}\wedge\dots\wedge e_{i_r}y(\lambda'))=(-1)^re_i\wedge e_{i_1}\wedge\dots\wedge e_{i_r}\pt{j}(y(\lambda')),\\
&E_{m+j,i}\cdot (e_{i_1}\wedge\dots\wedge e_{i_r}y(\lambda'))\\
&=\left\{
\begin{aligned}
&0 & \text{if} &\ i \notin \{i_1,\dots,i_r\},\\
&(-1)^{r-s}e_{i_1}\wedge\dots\wedge\hat e_{i_s}\wedge\dots\wedge e_{i_r}y(\lambda'+e_j) & \text{if}&\ i=i_s\in\{i_1,\dots,i_r\},
\end{aligned}\right.
\end{split}\end{equation*}
for any $i\in\{1,\dots,m\},j\in\{1,\dots,n\}$ and $e_{i_1}\wedge\dots\wedge e_{i_r}y(\lambda')\in M(\lambda)\setminus\{0\}.$
\end{lemma}

Consider the following space
\begin{eqnarray*}
&\span\{e_{i_1}\wedge\dots\wedge e_{i_r}y(\lambda')\mid r\in\{0,\dots,m\},\ i_1,\dots,i_r\in\{1,\dots,m\},\\
&m+|\lambda|=r+|\lambda'|,\lambda'\in\supp(W(\lambda))\}.
\end{eqnarray*}
It is easy to see that the space is a $\Z$-graded simple $\gl(m,n)$-submodule of $M(\lambda)$ with the top non-zero graded component being the $\gl(m,n)_0$-module $V(m)\otimes N(\lambda)$. Hence, the space is just $L(V(m)\otimes N(\lambda))$.
Similarly, for any $r\in\{0,\dots,m\}$, we have
\begin{eqnarray*}
&L(V(r)\otimes\C)=\span\{e_{i_1}\wedge\dots\wedge e_{i_{r'}}y(\alpha)\mid r'\in\{0,\dots, r\},\ i_1,\dots,i_{r'}\in\{1,\dots,m\},\\
&r=r'+|\alpha|,\alpha\in\supp(W({\bf 0}))\}.
\end{eqnarray*}
So for $r$ and $\lambda$ satisfying that $\lambda={\bf 0}$ if $r\neq m$, $L(V(r)\otimes N(\lambda))$ is a $\gl(m,n)$-submodule of $M(\lambda)$.

For $i\in\{1,\dots,m\}$ and $j\in\{1,\dots,n\}$, define the actions $e_i\wedge$ and $t_j\wedge$ on $M(\lambda)$ by:
\begin{eqnarray*}
&e_i\wedge(e_{i_1}\wedge\dots\wedge e_{i_r}y(\lambda'))=e_i\wedge e_{i_1}\wedge\dots\wedge e_{i_r}y(\lambda'),\\
&t_j\wedge(e_{i_1}\wedge\dots\wedge e_{i_r}y(\lambda'))=(-1)^re_{i_1}\wedge\dots\wedge e_{i_r}y(\lambda'+e_j),\forall \lambda'\in\supp(W(\lambda)).
\end{eqnarray*}
Let $P$ be a simple module over $K_{m,n}^+$ (resp. $K_{m,n}$). Define $\sigma_\lambda:F(P,M(\lambda))\rightarrow F(P,M(\lambda))$ by
\begin{equation}\label{defsig}
u\otimes v\mapsto\sum_{s=1}^m\pt{s}u\otimes e_s\wedge v+(-1)^{|u|-1}\sum_{l=1}^n\pxi{l}u\otimes t_l\wedge v,\forall u\in P,v\in M(\lambda).
\end{equation}
Then we have the following two lemmas, where the detailed proof of Lemma \ref{4.3} can be seen in Appendix.

\begin{lemma}\label{siglam}
\begin{itemize}
\item[(1)]$\sigma_\lambda^2=0.$
\item[(2)]$\sigma_\lambda(F(P,L(V(r)\otimes N(\lambda))))=0$ if and only if $(r;\lambda)=(m; (-1,\dots,-1))$.
\end{itemize}
\end{lemma}
\begin{proof}
(1). Let $u\in P$ and $v\in M(\lambda)$, then
\begin{eqnarray*}
&&\sigma_\lambda^2(u\otimes v)\\
&=&\sum_{s'=1}^m\sum_{s=1}^m\pt{s'}\pt{s}u\otimes e_{s'}\wedge(e_s\wedge v)
+(-1)^{|u|-1}\sum_{l'=1}^n\sum_{s=1}^m\pxi{l'}\pt{s}u\otimes t_{l'}\wedge(e_s\wedge v)\\
&&+(-1)^{|u|-1}\sum_{s'=1}^m\sum_{l=1}^n\pt{s'}\pxi{l}u\otimes e_{s'}\wedge(t_l\wedge v)
-\sum_{l'=1}^n\sum_{l=1}^n\pxi{l'}\pxi{l}u\otimes t_{l'}\wedge(t_l\wedge v)\\
&=&0.
\end{eqnarray*}
Note that the second equation holds because $e_{s'}\wedge(e_s\wedge v)=-e_s\wedge(e_{s'}\wedge v),t_{l'}\wedge(e_s\wedge v)=-e_s\wedge(t_{l'}\wedge v)$ and $t_{l'}\wedge(t_l\wedge v)=t_l\wedge(t_{l'}\wedge v)$.

(2). For any $i\in\{1,\dots,m\}$, we have
$$u=\pt{i}(t_iu)-t_i(\pt{i}(u)),\forall u\in P.$$
Then $\pt{i}(P)\neq 0$. Otherwise $P=0$, which is a contradiction. Similarly, for any $j\in\{1,\dots,n\}$, we have $\pxi{j}(P)\neq 0$.

Obviously, $e_i\wedge L(V(r)\otimes N(\lambda))=0$ for all $i\in\{1,\dots,m\}$ and $t_j\wedge L(V(r)\otimes N(\lambda))=0$ for all $j\in\{1,\dots,n\}$ if and only if $(r;\lambda)=(m; (-1,\dots,-1))$. So the conclusion holds.
\end{proof}

\begin{lemma}\label{4.3}
Let $P$ be a simple module over $K_{m,n}^+$ (resp. $K_{m,n}$) and $\lambda\in\C^n$. Then the map $\sigma_\lambda:F(P,M(\lambda))\rightarrow F(P,M(\lambda))$ is a homomorphism of $W_{m,n}^+$-modules (resp. $W_{m,n}$-modules).
\end{lemma}
Let $F(P,r,\lambda)=\sigma_\lambda(F(P,M(\lambda)))\cap F(P,L(V(r)\otimes N(\lambda)))$. Then $F(P,r,\lambda)$ is a submodule of $F(P,L(V(r)\otimes N(\lambda)))$. Clearly, $F(P,0,{\bf 0})=0$.
\begin{lemma}
Let $P$ be a simple module over $K_{m,n}^+$ (resp. $K_{m,n}$), $r\in\{0,\dots,m\}$ and $\lambda\in\C^n$ such that $\lambda={\bf 0}$ if $r\neq m$. Then $F(P,L(V(r)\otimes N(\lambda)))$ is not simple for all $r,\lambda$ except $(r;\lambda)=(0;{\bf 0})$ and $(r;\lambda)=(m;(-1,\dots,-1))$.
\end{lemma}
\begin{proof}
By Lemma \ref{siglam}(1), $\sigma_\lambda(F(P,r,\lambda))=0$ for all $(r;\lambda)$. Then from Lemma \ref{siglam}(2), $F(P,L(V(r)\otimes N(\lambda)))\neq F(P,r,\lambda)$ for $(r;\lambda)\neq(m;(-1,\dots,-1))$. Also, by the definitions of $F(P,r,\lambda)$ and $\sigma_\lambda$, $F(P,r,\lambda)\neq 0$ for $(r;\lambda)\neq(0;{\bf 0})$. So $F(P,r,\lambda)$ is a nonzero proper submodule of $F(P,L(V(r)\otimes N(\lambda)))$ for $(r;\lambda)\neq(0;{\bf 0})$ and $(m;(-1,\dots,-1))$.
\end{proof}
\begin{lemma}\label{mod}
Let $u\in P$ and $v\in L(V(r)\otimes N(\lambda))$, then we have
\begin{equation*}\begin{split}
t^\alpha\xi_I\pt{i}\cdot( u\otimes v)\equiv&\sum_{s=1}^mt^\alpha\xi_I\pt{s}u\otimes(\delta_{s,i}-E_{s,i})v\\
&+(-1)^{|u|}\sum_{l=1}^nt^\alpha\xi_I\pxi{l}u\otimes E_{m+l,i}v\ \text{mod}\ F(P,r,\lambda),\\
t^\alpha\xi_I\pxi{j}\cdot (u\otimes v)\equiv&(-1)^{|u|+1}\sum_{s=1}^mt^\alpha\xi_I\pt{s}u\otimes E_{s,m+j}v\\
&+\sum_{l=1}^nt^\alpha\xi_I\pxi{l}u\otimes(\delta_{l,j}+E_{m+l,m+j})v\ \text{mod}\ F(P,r,\lambda),
\end{split}\end{equation*}
for any $i\in\{1,\dots,m\}$ and $j\in\{1,\dots,n\}$.
\end{lemma}
\begin{proof}
Let $u'\in P$ and $v'=e_{i_1}\wedge\dots\wedge e_{i_{r'}}y(\lambda')\in L(V(r)\otimes N(\lambda))\setminus\{0\}$. We first claim that
\begin{eqnarray*}
a:&=&\sum_{s=1}^m\pt{s}u'\otimes E_{s,i}v'+(-1)^{|u'|-1}\sum_{l=1}^n\pxi{l}u'\otimes E_{m+l,i}v'\in F(P,r,\lambda),\\
b:&=&\sum_{s=1}^m\pt{s}u'\otimes E_{s,m+j}v'+(-1)^{|u'|-1}\sum_{l=1}^n\pxi{l}u'\otimes E_{m+l,m+j}v'\in F(P,r,\lambda).
\end{eqnarray*}

If $i\notin\{i_1,\dots,i_{r'}\}$, then obviously $a=0$ because $E_{s,i}v'=0$ and $E_{m+l,i}v'=0$.

If $i\in\{i_1,\dots,i_{r'}\}$, without loss of generality, we assume that $i=i_1$. Then
\begin{eqnarray*}
a&=&\sum_{s=1}^m\pt{s}u'\otimes e_s\wedge(e_{i_2}\wedge\dots\wedge e_{i_{r'}}y(\lambda'))\\
&&+(-1)^{|u'|-1}\sum_{l=1}^n\pxi{l}u'\otimes t_l\wedge(e_{i_2}\wedge\dots\wedge e_{i_{r'}}y(\lambda'))\in F(P,r,\lambda),\\
b&=&\lambda'_j(-1)^{r'}\sum_{s=1}^m\pt{s}u'\otimes e_s\wedge(e_{i_1}\wedge\dots\wedge e_{i_{r'}}y(\lambda'-e_j))\\
&&+(-1)^{|u'|-1+r'}\lambda'_j\sum_{l=1}^n\pxi{l}u'\otimes t_l\wedge(e_{i_1}\wedge\dots\wedge e_{i_{r'}}y(\lambda'-e_j))\in F(P,r,\lambda).
\end{eqnarray*}
Thus, the claim stands correctly. Therefore,
\begin{equation*}\begin{split}
&t^\alpha\xi_I\pt{i}\cdot (u\otimes v)\\
=&\sum_{s=1}^mt^\alpha\xi_I\pt{s}u\otimes(\delta_{s,i}-E_{s,i})v+(-1)^{|u|}\sum_{l=1}^nt^\alpha\xi_I\pxi{l}u\otimes E_{m+l,i}v\\
&+\sum_{s=1}^m\pt{s}t^\alpha\xi_Iu\otimes E_{s,i}v+(-1)^{|I|+|u|-1}\sum_{l=1}^n\pxi{l}t^\alpha\xi_Iu\otimes E_{m+l,i}v\\
\equiv&\sum_{s=1}^mt^\alpha\xi_I\pt{s}u\otimes(\delta_{s,i}-E_{s,i})v+(-1)^{|u|}\sum_{l=1}^nt^\alpha\xi_I\pxi{l}u\otimes E_{m+l,i}v\ \text{mod}\ F(P,r,\lambda),\\
&t^\alpha\xi_I\pxi{j}\cdot (u\otimes v)\\
=&(-1)^{|u|+1}\sum_{s=1}^mt^\alpha\xi_I\pt{s}u\otimes E_{s,m+j}v+\sum_{l=1}^nt^\alpha\xi_I\pxi{l}u\otimes(\delta_{l,j}+E_{m+l,m+j})v\\
\end{split}\end{equation*}
\begin{equation*}\begin{split}
&+(-1)^{|u|}\sum_{s=1}^m\pt{s}t^\alpha\xi_Iu\otimes E_{s,m+j}v+(-1)^{|I|-1}\sum_{l=1}^n\pxi{l}t^\alpha\xi_Iu\otimes E_{m+l,m+j}v\\
\equiv&(-1)^{|u|+1}\sum_{s=1}^mt^\alpha\xi_I\pt{s}u\otimes E_{s,m+j}v+\sum_{l=1}^nt^\alpha\xi_I\pxi{l}u\otimes(\delta_{l,j}+E_{m+l,m+j})v\ \text{mod}\ F(P,r,\lambda).
\end{split}\end{equation*}
\end{proof}

Let $\widetilde{F}(P,r,\lambda)=\{w\in F(P,L(V(r)\otimes N(\lambda)))\mid W_{m,n}^+w\subset F(P,r,\lambda)\}$ for $W_{m,n}^+$ (resp. $\widetilde{F}(P,r,\lambda)=\{w\in F(P,L(V(r)\otimes N(\lambda)))\mid W_{m,n}w\subset F(P,r,\lambda)\}$ for $W_{m,n}$). It is clear that $\widetilde{F}(P,r,\lambda)$ is a submodule of $F(P,L(V(r)\otimes N(\lambda)))$ containing $F(P,r,\lambda)$. Moreover, $\widetilde{F}(P,r,\lambda)/F(P,r,\lambda)$ is either zero or a trivial module.

\begin{lemma}\label{tildeF}
$\widetilde{F}(P,r,\lambda)=\{w\in F(P,L(V(r)\otimes N(\lambda)))\mid W_{m,n}^+w\subset \widetilde{F}(P,r,\lambda)\}$ for $W_{m,n}^+$ (resp. $\widetilde{F}(P,r,\lambda)=\{w\in F(P,L(V(r)\otimes N(\lambda)))\mid W_{m,n}w\subset \widetilde{F}(P,r,\lambda)\}$ for $W_{m,n}$).
\end{lemma}
\begin{proof}
We will prove this only for $W_{m,n}^+$ and this proof is also valid for $W_{m,n}$.

Suppose $m\geqslant 1$ or $n\geqslant 2$, then $[W_{m,n}^+,W_{m,n}^+]=W_{m,n}^+$ since $W_{m,n}^+$ is simple. For any $w\in F(P,L(V(r)\otimes N(\lambda)))$ with $W_{m,n}^+w\subset \widetilde{F}(P,r,\lambda)$, there is $W_{m,n}^+w\subset F(P,r,\lambda)$. So the conclusion holds for $m\geqslant 1$ or $n\geqslant 2$.

Suppose $m=0$ and $n=1$, then $P=\Lambda(1)$ and $N(\lambda)=\C v$ with $E_{1,1} v=\lambda_1v$. The action of $W_{0,1}^+$ on $F(\Lambda(1),\C v)$ is given by
\begin{equation*}
\pxi{1}\cdot (1\otimes v)=0,\ \xi_1\pxi{1}\cdot (1\otimes v)=\lambda_1\otimes v,\ \pxi{1}\cdot (\xi_1\otimes v)=1\otimes v,\ \xi_1\pxi{1}\cdot (\xi_1\otimes v)=(1+\lambda_1)\xi_1\otimes v.
\end{equation*}
\begin{itemize}
\item[] If $\lambda_1=0$, then $F(\Lambda(1),r,\lambda)=0$ and $\widetilde{F}(\Lambda(1),r,\lambda)=\span\{1\otimes v\}$.
\item[] If $\lambda_1=-1$, then $F(\Lambda(1),r,\lambda)=\span\{1\otimes v\}$ and $\widetilde{F}(\Lambda(1),r,\lambda)=F(\Lambda(1),\C v)$.
\item[] If $\lambda_1\neq 0,-1$, then $F(\Lambda(1),r,\lambda)=\widetilde{F}(\Lambda(1),r,\lambda)=\span\{1\otimes v\}$.
\end{itemize}
In a word, there is $\widetilde{F}(\Lambda(1),r,\lambda)=\{w\in F(\Lambda(1),\C v)\mid W_{0,1}^+w\subset \widetilde{F}(\Lambda(1),r,\lambda)\}.$
\end{proof}

\begin{lemma}\label{F/tildeF}
Let $P$ be a simple module over $K_{m,n}^+$ (resp. $K_{m,n}$), $r\in\{0,\dots,m\}$ and $\lambda\in\C^n$ such that $\lambda={\bf 0}$ if $r\neq m$. Then the $W_{m,n}^+$-module (resp. $W_{m,n}$-module) $F(P,L(V(r)\otimes N(\lambda)))/\widetilde{F}(P,r,\lambda)$ is either zero or simple.
\end{lemma}
\begin{proof}
We will prove this only for $W_{m,n}^+$ and this proof is also valid for $W_{m,n}$. Suppose that $\widetilde{F}(P,r,\lambda)\neq F(P,L(V(r)\otimes N(\lambda)))$. Let $F'$ be a submodule of $F(P,L(V(r)\otimes N(\lambda)))$ with $\widetilde{F}(P,r,\lambda)\subsetneq F'$.

Let $w=\sum_{p=1}^qu_p\otimes v_p\in F'\setminus\widetilde{F}(P,r,\lambda)$. By Lemma \ref{tildeF}, there is a $t^\beta\xi_I\partial\in W_{m,n}^+$ with $\partial\in\{\pt{1},\dots,\pt{m},\pxi{1},\dots,\pxi{n}\}$ such that $t^\beta\xi_I\partial\cdot w\notin \widetilde{F}(P,r,\lambda)$. By Lemma \ref{mod} and $F(P,r,\lambda)\subseteq \widetilde{F}(P,r,\lambda)$, there are $v_{p,p'},p\in\{1,\dots,q\},p'\in\{1,\dots,m+n\}$ such that for any $\alpha\in\Z_+^m$ and $I\subset J\subset\{1,\dots,n\}$, we have
$$t^{\alpha+\beta}\xi_J\partial\cdot w\equiv\sum_{p=1}^q\sum_{s=1}^mt^{\alpha+\beta}\xi_J\pt{s}u_p\otimes v_{p,s}+\sum_{p=1}^q\sum_{l=1}^nt^{\alpha+\beta}\xi_J\pxi{l}u_p\otimes v_{p,m+l}\ \text{mod}\ \widetilde{F}(P,r,\lambda).$$
Since $t^{\alpha+\beta}\xi_J\partial\cdot w\in F'$ and $\widetilde{F}(P,r,\lambda)\subsetneq F'$, we get
$\sum_{p=1}^q\sum_{s=1}^mt^{\alpha+\beta}\xi_J\pt{s}u_p\otimes v_{p,s}+\sum_{p=1}^q\sum_{l=1}^nt^{\alpha+\beta}\xi_J\pxi{l}u_p\otimes v_{p,m+l}\in F'$. It can be further obtained that
$$\sum_{p=1}^q\sum_{s=1}^mft^\beta\xi_I\pt{s}u_p\otimes v_{p,s}+\sum_{p=1}^q\sum_{l=1}^nft^\beta\xi_I\pxi{l}u_p\otimes v_{p,m+l}\in F',\forall f\in A_{m,n}^+.$$
Applying $\pt{i},i\in\{1,\dots,m\}$ and $\pxi{j},j\in\{1,\dots,n\}$ to the above element repeatedly, we have
$$\sum_{p=1}^q\sum_{s=1}^mxt^\beta\xi_I\pt{s}u_p\otimes v_{p,s}+\sum_{p=1}^q\sum_{l=1}^nxt^\beta\xi_I\pxi{l}u_p\otimes v_{p,m+l}\in F', \forall x\in K_{m,n}^+.$$
Note that
$$0\neq t^\beta\xi_I\partial\cdot w\equiv\sum_{p=1}^q\sum_{s=1}^mt^\beta\xi_I\pt{s}u_p\otimes v_{p,s}+\sum_{p=1}^q\sum_{l=1}^nt^\beta\xi_I\pxi{l}u_p\otimes v_{p,m+l}\ \text{mod}\ \widetilde{F}(P,r,\lambda).$$
By the Jacobson density theorem, we have $P\otimes v\subset F'$ for some $0\neq v\in L(V(r)\otimes N(\lambda))$. Since the space $M'=\{v\in L(V(r)\otimes N(\lambda))\mid P\otimes v\subset F'\}$ is a $\gl(m,n)$-submodule of $L(V(r)\otimes N(\lambda))$, $M'=L(V(r)\otimes N(\lambda))$, which implies $F'=F(P,L(V(r)\otimes N(\lambda)))$.
\end{proof}
\begin{theorem}\label{subquotients}
Let $P$ be a simple module over $K_{m,n}^+$ (resp. $K_{m,n}$), $r\in\{0,\dots,m\}$ and $\lambda\in\C^n$ such that $\lambda={\bf 0}$ if $r\neq m$.
\begin{itemize}
\item[(1)]$\widetilde{F}(P,r,\lambda)=\ker\sigma_\lambda\cap F(P,L(V(r)\otimes N(\lambda)))$ if $r\neq m$ or $\lambda\neq(-1,\dots,-1)$;
\item[(2)]$F(P,L(V(r)\otimes N(\lambda)))/\widetilde{F}(P,r,\lambda)$ is simple if $r\neq m$ or $\lambda\neq(-1,\dots,-1)$;
\item[(3)]$F(P,r,\lambda)$ is simple if $r\neq 0$ or $\lambda\neq {\bf 0}$;
\item[(4)]$F(P,L(V(0)\otimes N({\bf 0})))=P$ is simple if and only if $P$ is not isomorphic to the natural module $A_{m,n}^+$ (resp. $A_{m,n}$) up to a parity-change. Moreover,  the $K_{m,n}^+$-module $A_{m,n}^+/\C$ (resp. $A_{m,n}/\C$) is simple;
\item[(5)]Suppose $\lambda=(-1,\dots,-1)$, then $F(P,L(V(m)\otimes N(\lambda)))$ is simple if and only if $P=\sum_{s=1}^m\pt{s}P+\sum_{l=1}^n\pxi{l}P$. For $P\neq\sum_{s=1}^m\pt{s}P+\sum_{l=1}^n\pxi{l}P$, $F(P,L(V(m)\otimes N(\lambda)))/F(P,m,\lambda)$ is a trivial module.
\end{itemize}
\end{theorem}
\begin{proof}
We will prove this results only for $W_{m,n}^+$ since the proof is valid also for $W_{m,n}$.

(1)-(2). Suppose that there is a $w\in\widetilde{F}(P,r,\lambda)$ such that $\sigma_\lambda(w)\neq 0$. For any $x\in W_{m,n}^+$, $xw\in F(P,r,\lambda)$. Then $x\sigma_\lambda(w)=\sigma_\lambda(xw)=0$ by Lemma \ref{siglam}(1). Let $\sigma_\lambda(w)=\sum_{p=1}^qu_p\otimes v_p$ with $v_1,\dots,v_q$ linearly independent, then
$$0=\pt{i}\cdot\sigma_\lambda(w)=\sum_{p=1}^q\pt{i}u_p\otimes v_p,\ 0=\pxi{j}\cdot\sigma_\lambda(w)=\sum_{p=1}^q\pxi{j}u_p\otimes v_p.$$
It follows that $\pt{i}u_p=0,\ \pxi{j}u_p=0$ for any $i\in\{1,\dots,m\},j\in\{1,\dots,n\}$ and $p\in\{1,\dots,q\}$.
Since $P$ is simple module over $K_{m,n}^+$, we have $P=K_{m,n}^+u_1$. Also, it is easy to see
$K_{m,n}^+/(\sum\limits_{i=1}^m K_{m,n}^+\pt{i}+\sum\limits_{j=1}^n K_{m,n}^+\pxi{j})\cong A_{m,n}^+$. Therefore, we get $P\cong A_{m,n}^+$ or $\Pi(A_{m,n}^+)$ and $u_p\in\C$ for $p\in \{1,\dots,q\}$. So $q=1$ and $\sigma_\lambda(w)=1\otimes v$ for some $v\in M(\lambda)$. Then for any $i\in\{1,\dots,m\}$,
$0=t_i\pt{i}\cdot \sigma_\lambda(w)=1\otimes E_{i,i}v,$
which implies $E_{i,i}v=0$. Similarly, $E_{m+j,m+j}v=0$ for $j\in\{1,\dots,n\}$. Therefore, $1\otimes v\in F(P,L(V(0)\otimes \C))$, which contradicts with the fact that $F(P,0,{\bf 0})=0$. So, $\widetilde{F}(P,r,\lambda)\subset\ker\sigma_\lambda$.

By Lemma \ref{siglam}, $F(P,L(V(r)\otimes N(\lambda)))\neq \ker\sigma_\lambda\cap F(P,L(V(r)\otimes N(\lambda)))$ if $r\neq m$ or $\lambda\neq(-1,\dots,-1)$. Then (1) and (2) follow from Lemma \ref{F/tildeF}.

(3). If $r\neq 0$ and $\lambda={\bf 0}$, then $F(P,r,{\bf 0})=\sigma_0(F(P,L(V(r-1)\otimes\C)))\neq 0$. By (1) and (2),
$$F(P,r,{\bf 0})\cong F(P,L(V(r-1)\otimes\C))/\widetilde{F}(P,r-1,{\bf 0})$$
is simple. If $r=m$ and $\lambda\neq {\bf 0}$, assume that $\lambda_1\neq 0$. Then $F(P,m,\lambda)=\sigma_\lambda(F(P,L(V(m)\otimes N(\lambda-e_1))))\neq 0$. By (1) and (2),
$$F(P,m,\lambda)\cong F(P,L(V(m)\otimes N(\lambda-e_1)))/\widetilde{F}(P,m,\lambda-e_1)$$
is simple. Hence we get (3).

(4). By (\ref{P}), $P\cong P_1\otimes\Lambda(n)$, where $P_1$ is a simple module over $K_{m,0}^+$. By Theorem 4 in \cite{TZ}, $P_1$ is a simple module over $W_{m,0}^+$ if and only if $P_1$ is not isomorphic to the natural $K_{m,0}^+$-module $A_{m,0}^+$ up to a parity-change.

If $P_1$ is isomorphic to $A_{m,0}^+$ or $\Pi(A_{m,0}^+)$, then $P$ is isomorphic to $A_{m,n}^+$ or $\Pi(A_{m,n}^+)$. Clearly, the $W_{m,n}^+$-module $A_{m,n}^+$ is not simple but $A_{m,n}^+/\C$ is.

If $P_1$ is not isomorphic to $A_{m,0}^+$ or $\Pi(A_{m,0}^+)$, then $P_1$ is a simple module over $W_{m,0}^+$. Let $Q$ be a nonzero $W_{m,n}^+$-submodule of $P$. By the Jacobson density theorem, $P_1\otimes u\subset Q$ for some nonzero $u\in \Lambda(n)$. Let $u'\in P_1$ with $\pt{1}u'\neq 0$. Applying $\pxi{j}, j\in\{1,\dots,n\}$ to $u'\otimes u$, we get $u'\otimes 1\in Q$. Then
$$\xi_{I_n}\pt{1}\cdot (u'\otimes 1)=\pt{1}u'\otimes \xi_{I_n}\in Q,$$
where $I_n=\{1,\dots,n\}$. Applying elements in $W_{m,0}^+$ and $\pxi{j}, j\in\{1,\dots,n\}$ to $\pt{1}u'\otimes \xi_{I_n}$, we get $P_1\otimes \Lambda(n)\subset Q$. Thus $P$ is a simple module over $W_{m,n}^+$ and (4) is proved.

(5). Since $\lambda=(-1,\dots, -1)$, $L(V(m)\otimes N(\lambda))=\C v$ with $v=e_1\wedge\dots\wedge e_my(\lambda)$.
From (\ref{defsig}) we know that
$$F(P,m,\lambda)=(\sum_{s=1}^m\pt{s}P+\sum_{l=1}^n\pxi{l}P)\otimes v.$$
By (3), $F(P,m,\lambda)$ is a simple submodule of $F(P,L(V(m)\otimes N(\lambda)))$. So $F(P,L(V(m)\otimes N(\lambda)))$ is simple if and only if $P=\sum_{s=1}^m\pt{s}P+\sum_{l=1}^n\pxi{l}P$.
For any $u\in P$,
$$t^\alpha\xi_I\pt{i}\cdot (u\otimes v)=t^\alpha\xi_I\pt{i}u\otimes v+\pt{i}(t^\alpha)\xi_Iu\otimes v=\pt{i}t^\alpha\xi_Iu\otimes v\in F(P,m,\lambda),$$
$$t^\alpha\xi_I\pxi{j}\cdot (u\otimes v)=t^\alpha\xi_I\pxi{j}u\otimes v+(-1)^{|I|}t^\alpha\pxi{j}(\xi_I)u\otimes v=(-1)^{|I|}\pxi{j}t^\alpha\xi_Iu\otimes v\in F(P,m,\lambda).$$
Therefore, $F(P,L(V(m)\otimes N(\lambda)))/F(P,m,\lambda)$ is a trivial module if $P\neq\sum_{s=1}^m\pt{s}P+\sum_{l=1}^n\pxi{l}P$.
\end{proof}

Finally, we will give the isomorphism criterion for two irreducible modules $F(P, M)$. Let $M$ be a module over a Lie superalgebra or an associative superalgebra $\g$. Then we can get a new module, denoted by $M^T$, over $\g$ by $xv=(-1)^{|x|}xv, \forall x\in \g, v\in M^T$.
\begin{theorem}
Let $P, P'$ be simple $K_{m,n}^+$-module (resp. $K_{m,n}$-module), and $M, M'$ be simple $\gl(m, n)$-module. Suppose that $M$ is not isomorphic to
$L(V_1\otimes V_2)$, where $V_1=V(r)$ for $r\in\{0, \cdots, m\}$ and $V_2$ is a simple $\gl_n$-module. Then $F(P, M)\cong F(P', M')$ if and only
if $P\cong P'$ and $M\cong M'$ (or $P\cong \Pi(P')$ and $M\cong \Pi(M'^T)$).
\end{theorem}
\begin{proof}
We will prove this result only for $W_{m,n}^+$ since the proof is valid also for $W_{m,n}$.
The sufficiency is obvious. Next we prove the necessary. Assume that the following elements in $P, M$ and $K_{m,n}^+$ are all homogeneous. Suppose that
$$\psi: F(P, M)\rightarrow F(P', M')$$
is an isomorphism of $W_{m,n}^+$-modules. Let $u\otimes v$ be a nonzero element in $F(P, M)$ and assume that
$$\psi(u\otimes v)=\sum_{p=1}^q u'_p\otimes v'_p,$$
where $u'_1, \cdots, u'_q$ are linearly independent. As in the proof of the Claim 1 in Lemma \ref{V_1}, we have
\begin{eqnarray}\label{3.4.1}
\psi(xu\otimes(\delta_{r,i}E_{r,k}-E_{r,i}E_{r,k})v)=
\sum_{p=1}^q xu'_p\otimes(\delta_{r,i}E_{r,k}-E_{r,i}E_{r,k})v'_p
\end{eqnarray}
for any $i, k, r\in\{1,\dots,m\}$ and $x\in K_{m,n}^+$.
By rechoosing $u\otimes v$, since $M$ is not isomorphic to $L(V_1\otimes V_2)$ with $V_1=V(r)$ for $r\in\{0,\cdots, m\}$,
without loss of generality, we can assume that $(\delta_{r,i}E_{r,k}-E_{r,i}E_{r,k})v'_1 \neq 0$
for some $i, k, r$. Furthermore, since $u'_1, \cdots, u'_q$ are linearly independent, from the Jacobson density theorem,
we can find $y\in K_{m,n}^+$ so that $yu'_p=\delta_{1,p}u'_1, p\in\{1,\cdots, q\}$. Let $x=y$ in (\ref{3.4.1}), we have
\begin{eqnarray*}
\psi(yu\otimes(\delta_{r,i}E_{r,k}-E_{r,i}E_{r,k})v)=u'_1\otimes(\delta_{r,i}E_{r,k}-E_{r,i}E_{r,k})v'_1\neq 0,
\end{eqnarray*}
which implies that $yu\neq 0$ and $(\delta_{r,i}E_{r,k}-E_{r,i}E_{r,k})v\neq 0$.
Now replacing $x$ with $xy$ in (\ref{3.4.1}), then replacing $yu$ with $u$,
$(\delta_{r,i}E_{r,k}-E_{r,i}E_{r,k})v$ with $v$, $u'_1$ with $u'$,
and $(\delta_{r,i}E_{r,k}-E_{r,i}E_{r,k})v'_1$ with $v'$, we get
\begin{equation}\label{3.4.2}
\psi(xu\otimes v)=xu'\otimes v', \forall x\in K_{m,n}^+.
\end{equation}
Let $\mbox{Ann}_{K_{m,n}^+}(u)=\{x\in K_{m,n}^+\mid xu=0\}$. Since $\psi$ is an isomorphism, (\ref{3.4.2}) implies
that $\mbox{Ann}_{K_{m,n}^+}(u)=\mbox{Ann}_{K_{m,n}^+}(u')$. Also, we have
$$K_{m,n}^+u=P\ \text{and}\ K_{m,n}^+u'=P'$$
because $P$ and $P'$ are simple $K_{m,n}^+$-modules.
If $u\in P_{\bar{0}}$ (resp. $u'\in P'_{\bar{0}}$), we have $K_{m,n}^+/{\mbox{Ann}_{K_{m,n}^+}(u)}\cong P$ (resp. $K_{m,n}^+/{\mbox{Ann}_{K_{m,n}^+}(u')}\cong P'$).
If $u\in P_{\bar{1}}$ (resp. $u'\in P'_{\bar{1}}$), we have $K_{m,n}^+/{\mbox{Ann}_{K_{m,n}^+}(u)}\cong \Pi(P)$ (resp. $K_{m,n}^+/{\mbox{Ann}_{K_{m,n}^+}(u')}\cong \Pi(P')$).
Therefore, it follows that $P\cong P' (\text{or}\ \Pi(P'))$ since $K_{m,n}^+/{\mbox{Ann}_{K_{m,n}^+}(u)}=K_{m,n}^+/{\mbox{Ann}_{K_{m,n}^+}(u')}$.

If $P\cong P'$, the map $\psi_1: P\rightarrow P'$ with $\psi_1(xu)=xu'$ gives the isomorphism of $K_{m,n}^+$-modules. Hence
\begin{equation}\label{3.4.3}
\psi(u\otimes v)=\psi_1(u)\otimes v', \forall u\in P.
\end{equation}
Now, by applying $t_i\pt{s}, \xi_j\pt{s}, t_i\pxi{l}$ and $\xi_j\pxi{l}$ to $\psi(u\otimes v)$ respectively, and combining with (\ref{3.4.3}), we deduce that
$$\psi(u\otimes E_{i,s}v)=\psi_1(u)\otimes E_{i,s}v',\ \psi(u\otimes E_{m+j,s}v)=\psi_1(u)\otimes E_{m+j,s}v',$$
$$\psi(u\otimes E_{i,m+l}v)=\psi_1(u)\otimes E_{i,m+l}v',\ \psi(u\otimes E_{m+j,m+l}v)=\psi_1(u)\otimes E_{m+j,m+l}v',$$
for any $i, s\in \{1, \cdots, m\}, j,l \in\{1, \cdots, n\}$ and $u\in P$. Consequently,
$$\psi(u\otimes yv)=\psi_1(u)\otimes yv', \forall u\in P, y\in U(\gl(m,n)).$$
This implies that $\mbox{Ann}_{U(\gl(m,n))}(v)=\mbox{Ann}_{U(\gl(m,n))}(v')$. Since $M$ and $M'$ are simple $U(\gl(m,n))$-modules, and
$$U(\gl(m,n))/{\mbox{Ann}_{U(\gl(m,n))}(v)}=U(\gl(m,n))/{\mbox{Ann}_{U(\gl(m,n))}(v')},$$
it follows that $M\cong M'$ or $M\cong \Pi(M')$. But since $|v|=|v'|$ in (\ref{3.4.3}), the case of $M\cong \Pi(M')$ can be ruled out.

If $P\cong \Pi(P')$, the map $\psi_2: P\rightarrow \Pi(P')$ with $\psi_2(xu)=xu'$ gives the isomorphism $K_{m,n}^+$-modules. Hence
\begin{equation}\label{3.4.4}
\psi(u\otimes v)=\psi_2(u)\otimes v', \forall u\in P.
\end{equation}
Note that $|v|=|v'|-1$ in (\ref{3.4.4}) since the parity of the vector $\psi_2(u)$ in $P'$ is $|u|+1$.
Similarly, by applying $t_i\pt{s}, \xi_j\pt{s}, t_i\pxi{l}$ and $\xi_j\pxi{l}$ to $\psi(u\otimes v)$ respectively, and combining with (\ref{3.4.4}), we get
$$\psi(u\otimes yv)=\psi_2(u)\otimes (-1)^{|y|}yv', \forall u\in P, y\in U(\gl(m,n)).$$
This also implies that $\mbox{Ann}_{U(\gl(m,n))}(v)=\mbox{Ann}_{U(\gl(m,n))}(v')$. So at this point, we can get $M\cong \Pi(M'^T)$.
\end{proof}

\section*{Appendix}
\setcounter{equation}{0}
\renewcommand{\theequation}{A\arabic{equation}}

\noindent{\bf Proof of Lemma \ref{4.1}.}
Let $v=e_{i_1}\wedge\dots\wedge e_{i_r}y(\lambda')\in M(\lambda)\setminus\{0\}$. To prove that $M(\lambda)$ is a $\gl(m,n)$-module, what we only need is to verify that the following seven equations hold since $M(\lambda)$ is already a $\gl(m,n)_0$-module. For any $i,i',k\in\{1,\dots,m\}$ and $j,j',l\in\{1,\dots,n\}$,
\begin{eqnarray}
(E_{i,k}E_{m+j,i'}-E_{m+j,i'}E_{i,k}-[E_{i,k},E_{m+j,i'}])\cdot v & = & 0,\label{A1}\\
(E_{i,k}E_{i',m+j}-E_{i',m+j}E_{i,k}-[E_{i,k},E_{i',m+j}])\cdot v & = & 0,\label{A2}\\
(E_{i,m+j}E_{m+j',i'}+E_{m+j',i'}E_{i,m+j}-[E_{i,m+j},E_{m+j',i'}])\cdot v & = & 0,\\
(E_{i,m+j}E_{i',m+j'}+E_{i',m+j'}E_{i,m+j}-[E_{i,m+j},E_{i',m+j'}])\cdot v & = & 0,\\
(E_{m+j,i}E_{m+j',i'}+E_{m+j',i'}E_{m+j,i}-[E_{m+j,i},E_{m+j',i'}])\cdot v & = & 0,\\
(E_{m+l,m+j'}E_{m+j,i'}-E_{m+j,i'}E_{m+l,m+j'}-[E_{m+l,m+j'},E_{m+j,i'}])\cdot v & = & 0,\\
(E_{m+l,m+j'}E_{i',m+j}-E_{i',m+j}E_{m+l,m+j'}-[E_{m+l,m+j'},E_{i',m+j}])\cdot v & = & 0.\label{A7}
\end{eqnarray}
Equtions (\ref{A1})-(\ref{A7}) can be directly verified by using the action of $\gl(m,n)_{\pm 1}$ on $M(\lambda)$. Here we only give a detailed proof of (\ref{A1}).

For any $i,i',k\in\{1,\dots,m\}$ and $j\in\{1,\dots,n\}$, we have
\begin{equation*}\begin{split}
&E_{i,k}E_{m+j,i'}\cdot (e_{i_1}\wedge\dots\wedge e_{i_r} y(\lambda'))\\
=&\left\{\begin{aligned}
&0&\text{if}&\ i'\notin\{i_1,\dots,i_r\},\\
&(-1)^{r-s'}E_{i,k}\cdot (e_{i_1}\wedge\dots\wedge \hat e_{i_{s'}}\wedge\dots\wedge e_{i_r}y(\lambda'+e_{j}))&\text{if}&\ i'=i_{s'}\in\{i_1,\dots,i_r\},
\end{aligned}\right.\\
&E_{m+j,i'}E_{i,k}\cdot (e_{i_1}\wedge\dots\wedge e_{i_r}y(\lambda'))\\
=&\left\{\begin{aligned}
&0&\text{if}&\ k\notin\{i_1,\dots,i_r\},\\
&E_{m+j,i'}\cdot((-1)^{s-1}e_i\wedge e_{i_1}\wedge\dots\wedge \hat e_{i_s}\wedge\dots\wedge e_{i_r}y(\lambda'))&\text{if}&\ k=i_s\in\{i_1,\dots,i_r\},
\end{aligned}\right.\\
&[E_{i,k},E_{m+j,i'}]\cdot (e_{i_1}\wedge\dots\wedge e_{i_r} y(\lambda')\\
=&\left\{\begin{aligned}
&0&\text{if}&\ k\notin\{i_1,\dots,i_r\},\\
&(-1)^{r-s-1}\delta_{i,i'}e_{i_1}\wedge\dots\wedge\hat e_{i_s}\wedge\dots\wedge e_{i_r}y(\lambda'+e_{j})&\text{if}&\ k=i_s\in\{i_1,\dots,i_r\}.
\end{aligned}\right.
\end{split}\end{equation*}

If $k\notin\{i_1,\dots,i_r\}$, then (\ref{A1}) holds.

If $k=i_s\in\{i_1,\dots,i_r\}$, suppose $i'\notin\{i_1,\dots,i_r\}$, then
$$E_{m+j,i'}E_{i,k}\cdot (e_{i_1}\wedge\dots\wedge e_{i_r}y(\lambda'))=(-1)^{r-s}\delta_{i,i'}e_{i_1}\wedge\dots\wedge\hat e_{i_s}\dots\wedge e_{i_r}y(\lambda'+e_{j}).$$
So (\ref{A1}) holds. Suppose $i'=i_{s'}\in\{i_1,\dots,i_r\}$, then
\begin{equation*}\begin{split}
&E_{i,k}E_{m+j,i'}\cdot (e_{i_1}\wedge\dots\wedge e_{i_r} y(\lambda'))\\
&=\left\{\begin{aligned}
&(-1)^{r-s'-s-1}e_i\wedge e_{i_1}\wedge\dots\wedge\hat e_{i_s}\dots\wedge\hat e_{i_{s'}}\dots\wedge e_{i_r}y(\lambda'+e_j)&\text{if}&\ s<s',\\
&0&\text{if}&\ s=s',\\
&(-1)^{r-s'-s}e_i\wedge e_{i_1}\wedge\dots\wedge\hat e_{i_{s'}}\dots\wedge\hat e_{i_s}\dots\wedge e_{i_r}y(\lambda'+e_j)&\text{if}&\ s>s'.
\end{aligned}\right.
\end{split}\end{equation*}
At this point, if $i\in\{i_1,\dots,i_r\}\setminus\{i_s\}$, we have
\begin{equation*}\begin{split}
E_{i,k}E_{m+j,i'}\cdot (e_{i_1}\wedge\dots\wedge e_{i_r} y(\lambda'))
=(-1)^{r-s-1}\delta_{i,i'}e_{i_1}\wedge\dots\wedge\hat e_{i_s}\dots\wedge e_{i_r}y(\lambda'+e_j)
\end{split}\end{equation*}
and $E_{m+j,i'}E_{i,k}\cdot (e_{i_1}\wedge\dots\wedge e_{i_r} y(\lambda'))=0$. So (\ref{A1}) holds. If $i\notin\{i_1,\dots,i_r\}\setminus\{i_s\}$, we have
\begin{equation*}\begin{split}
&E_{m+j,i'}E_{i,k}\cdot (e_{i_1}\wedge\dots\wedge e_{i_r} y(\lambda'))\\
&=\left\{\begin{aligned}
&(-1)^{r-s'-s-1}e_i\wedge e_{i_1}\wedge\dots\wedge\hat e_{i_s}\dots\wedge\hat e_{i_{s'}}\dots\wedge e_{i_r}y(\lambda'+e_j)&\text{if}&\ s<s',\\
&(-1)^{r-s}\delta_{i,j'}e_{i_1}\wedge\dots\wedge\hat e_{i_s}\dots\wedge e_{i_r}y(\lambda'+e_j)&\text{if}&\ s=s',\\
&(-1)^{r-s'-s}e_i\wedge e_{i_1}\wedge\dots\wedge\hat e_{i_{s'}}\dots\wedge\hat e_{i_s}\dots\wedge e_{i_r}y(\lambda'+e_j)&\text{if}&\ s>s'.
\end{aligned}\right.
\end{split}\end{equation*}
Therefore, (\ref{A1}) can still be obtained case by case.

\noindent{\bf Proof of Lemma \ref{4.3}.}
Let $u\in P$ and $v=e_{i_1}\wedge\dots\wedge e_{i_r}y(\lambda')\in M(\lambda)\setminus\{0\}$.
For any $\alpha\in\Z_+^m$ (or $\Z^m$), $i\in\{1,\dots,m\}, j\in\{1,\dots,n\}$ and $I\subset\{1,\dots,n\}$, we need to prove the following two facts:
$$A=\sigma_\lambda(t^\alpha\xi_I\pt{i}\cdot (u\otimes v))-t^\alpha\xi_I\pt{i}\cdot \sigma_\lambda(u\otimes v)=0,$$
$$B=\sigma_\lambda(t^\alpha\xi_I\pxi{j}\cdot (u\otimes v))-t^\alpha\xi_I\pxi{j}\cdot \sigma_\lambda(u\otimes v)=0.$$
We only give the detailed proof for $A=0$ and the other similar case is omited here. For simplicity, we denote $(-1)^{|I|+|u|-1}$ by $(-1)^{a}$ in the following equations.
\begin{eqnarray*}
&&\sigma_\lambda(t^\alpha\xi_I\pt{i}\cdot (u\otimes v))\\
&=&\sigma_\lambda\big(t^\alpha\xi_I\pt{i}u\otimes v+\sum_{s=1}^m\pt{s}(t^\alpha)\xi_Iu\otimes E_{s,i}v+(-1)^a\sum_{l=1}^n\pxi{l}(\xi_I)t^\alpha u\otimes E_{m+l,i}v\big)\\
&=&\sum_{s'=1}^m\pt{s'}t^\alpha\xi_I\pt{i}u\otimes e_{s'}\wedge v+(-1)^a\sum_{l'=1}^n\pxi{l'}t^\alpha\xi_I\pt{i}u\otimes t_{l'}\wedge v\\
&&+\sum_{s,s'=1}^m\pt{s'}\pt{s}(t^\alpha)\xi_Iu\otimes e_{s'}\wedge E_{s,i}v
+(-1)^a\sum_{l'=1}^n\sum_{s=1}^m\pxi{l'}\pt{s}(t^\alpha)\xi_Iu\otimes t_{l'}\wedge E_{s,i}v\\
&&+(-1)^a\sum_{s'=1}^m\sum_{l=1}^n\pt{s'}\pxi{l}(\xi_I)t^\alpha u\otimes e_{s'}\wedge E_{m+l,i}v-\sum_{l,l'=1}^n\pxi{l'}\pxi{l}(\xi_I)t^\alpha u\otimes t_{l'}\wedge E_{m+l,i}v,\\
&&t^\alpha\xi_I\pt{i}\cdot\sigma_\lambda(u\otimes v)\\
&=&\sum_{s'=1}^mt^\alpha\xi_I\pt{i}\pt{s'}u\otimes e_{s'}\wedge u+(-1)^a\sum_{s'=1}^m\sum_{l=1}^n\pxi{l}(\xi_I)t^\alpha\pt{s'}u\otimes E_{m+l,i}(e_{s'}\wedge v)\\
&&+(-1)^{|u|-1}\sum_{l'=1}^nt^\alpha\xi_I\pt{i}\pxi{l'}u\otimes t_{l'}\wedge v+(-1)^{|u|-1}\sum_{l'=1}^n\sum_{s=1}^m\pt{s}(t^\alpha)\xi_I\pxi{l'}u\otimes E_{s,i}(t_{l'}\wedge v)\\
&&+\sum_{s,s'=1}^m\pt{s}(t^\alpha)\xi_I\pt{s'}u\otimes E_{s,i}(e_{s'}\wedge v)+(-1)^{|I|-1}\sum_{l,l'=1}^n\pxi{l}(\xi_I)t^\alpha\pxi{l'}u\otimes E_{m+l,i}(t_{l'}\wedge v)\\
&=&\sum_{s'=1}^mt^\alpha\xi_I\pt{i}\pt{s'}u\otimes e_{s'}\wedge v+\sum_{s=1}^m\pt{s}(t^\alpha)\xi_I\pt{i}u\otimes e_s\wedge v\\
&&+\sum_{s,s'=1}^m\pt{s}(t^\alpha)\xi_I\pt{s'}u\otimes e_{s'}\wedge E_{s,i}v+(-1)^a\sum_{s'=1}^m\sum_{l=1}^n\pxi{l}(\xi_I)t^\alpha\pt{s'}u\otimes E_{m+l,i}(e_{s'}\wedge v)\\
&&+(-1)^{|u|-1}\sum_{l'=1}^nt^\alpha\xi_I\pt{i}\pxi{l'}u\otimes t_{l'}\wedge v+(-1)^{|u|-1}\sum_{l'=1}^n\sum_{s=1}^m\pt{s}(t^\alpha)\xi_I\pxi{l'}u\otimes t_{l'}\wedge E_{s,i}v\\
&&+(-1)^{|I|}\sum_{l,l'=1}^n\pxi{l}(\xi_I)t^\alpha\pxi{l'}u\otimes t_{l'}\wedge E_{m+l,i}v.
\end{eqnarray*}
Note that the last equation follows from $E_{s,i}(e_{s'}\wedge v)=\delta_{s',i}e_s\wedge v+e_{s'}\wedge E_{s,i}v,E_{s,i}(t_{l'}\wedge v)=t_{l'}\wedge E_{s,i}v$ and $E_{m+l,i}(t_{l'}\wedge v)=-t_{l'}\wedge E_{m+l,i}v$. Then
\begin{eqnarray*}
A&=&\sum_{s'=1}^m\sum_{s=1}^m\pt{s'}(\pt{s}(t^\alpha))\xi_Iu\otimes e_{s'}\wedge E_{s,i}v+(-1)^a\sum_{l'=1}^nt^\alpha\pt{i}\pxi{l'}(\xi_I)u\otimes t_{l'}\wedge v\\
&+&(-1)^{|I|+|u|}\sum_{s'=1}^m\sum_{l=1}^n\pxi{l}(\xi_I)t^\alpha\pt{s'}u\otimes E_{m+l,i}(e_{s'}\wedge v)\\
&+&(-1)^a\sum_{l'=1}^n\sum_{s=1}^m\pt{s}(t^\alpha)\pxi{l'}(\xi_I)u\otimes t_{l'}\wedge E_{s,i}v-\sum_{l,l'=1}^n\pxi{l'}(\pxi{l}(\xi_I))t^\alpha u\otimes t_{l'}\wedge E_{m+l,i}v\\
&+&(-1)^a\sum_{s'=1}^m\sum_{l=1}^n\pt{s'}\pxi{l}(\xi_I)t^\alpha u\otimes e_{s'}\wedge E_{m+l,i}v.
\end{eqnarray*}
If $i\notin\{i_1,\dots,i_r\}$, then
\begin{eqnarray*}
A&=&(-1)^{|I|+|u|}\sum_{l=1}^n\pxi{l}(\xi_I)t^\alpha\pt{i}u\otimes(-1)^re_{i_1}\wedge\dots\wedge e_{i_r}y(\lambda'+e_l)\\
&+&(-1)^a\sum_{l'=1}^nt^\alpha\pt{i}\pxi{l'}(\xi_I)u\otimes(-1)^re_{i_1}\wedge\dots\wedge e_{i_r}y(\lambda'+e_{l'})=0.\\
\end{eqnarray*}
If $i=i_k\in\{i_1,\dots,i_r\}$, then
\begin{eqnarray*}
A&=&(-1)^{|I|+|u|}\sum_{s'=1}^m\sum_{l=1}^n\pxi{l}(\xi_I)t^\alpha\pt{s'}u\otimes(1-\delta_{i,s'})(-1)^{r-k}\\
&&e_{s'}\wedge e_{i_1}\wedge\dots\wedge\hat e_{i_k}\dots\wedge e_{i_r}y(\lambda'+e_l)\\
&+&(-1)^a\sum_{l'=1}^nt^\alpha\pt{i}\pxi{l'}(\xi_I)u\otimes(-1)^re_{i_1}\wedge\dots\wedge e_{i_r}y(\lambda'+e_{l'})\\
&+&(-1)^a\sum_{l'=1}^n\sum_{s=1}^m\pt{s}(t^\alpha)\pxi{l'}(\xi_I)u\otimes(-1)^{k-1+r}e_s\wedge e_{i_1}\wedge\dots\wedge\hat e_{i_k}\dots\wedge e_{i_r}y(\lambda'+e_{l'})\\
&+&(-1)^a\sum_{s'=1}^m\sum_{l=1}^n\pt{s'}\pxi{l}(\xi_I)t^\alpha u\otimes(-1)^{r-k}e_{s'}\wedge e_{i_1}\wedge\dots\wedge\hat e_{i_k}\dots\wedge e_{i_r}y(\lambda'+e_l)=0.
\end{eqnarray*}
Hence, $\sigma_\lambda$ is a homomorphism of $W_{m,n}^+$-modules (resp. $W_{m,n}$-modules).

\noindent {\bf Acknowledgement} The authors would like to thank the professor R. L\"{u} for formulating the problem and for his help in preparation of this paper. Y. Xue is partially supported by NSF of China (Grant 11771122, 11971440, 11801390). Y. Wang is supported by NSF of China (Grant 11871052).

\
Y. Xue: Department of Mathematics, Soochow University, Suzhou, P. R. China, Email: yhxue00@stu.suda.edu.cn.

\
Y. Wang: School of Mathematics, Tianjin University, Tianjin, P. R. China, Email: wangyan09@tju.edu.cn, Corresponding author.


\begin{thebibliography}{00}
\bibitem{BBL}  G. Benkart, D. Britten, F. Lemire, Modules with bounded weight multiplicities for simple Lie algebras, Math. Z, 225 (1997), no. 2, 333-353.
\bibitem{BF1}  Y. Billig, V. Futorny, Classification of irreducible representations of Lie algebra of vector fields on a torus, J. Reine Angew. Math., 720 (2016), 199-216.
\bibitem{BF2}  Y. Billig, V. Futorny,  Classification of simple cuspidal modules for solenoidal Lie algebras, Israel J. Math., 222 (2017), no. 1, 109-123.
\bibitem{BFIK}  Y. Billig, V. Futorny, K. Iohara, I. Kashuba, Classification of simple strong Harish-Chandra $W(m,n)$-modules, arXiv:2006.05618.
\bibitem{CG} A. Cavaness, D. Grantcharov, Bounded weight modules of the Lie algebra of vector fields on $\C^2$, J. Algebra Appl., 16(2017), no. 12, 1750236.
\bibitem{CM}  C. Chen, V. Mazorchuk, Simple supermodules over Lie superalgebras, Trans. Amer. Math. Soc., 374 (2021), no. 2, 899–921.
\bibitem{DMP} I. Dimitrov, O. Mathieu, I. Penkov, On the stucture of weight modules, Trans. Amer. Math. Soc., 352 (2000), 2857-2869.
\bibitem{E1}  S. Eswara Rao, Irreducible representations of the Lie algebra of the diffeomorphisms of a d-dimensional torus, J. Algebra, 182 (1996), no. 2, 401-421.
\bibitem{E2}  S. Eswara Rao, Partial classification of modules for Lie algebra of diffeomorphisms of d-dimensional torus, J. Math. Phys., 45 (2004), 3322-3333.
\bibitem{GS}  D. Grantcharov, V. Serganova, Category of ${\rm sp}(2n)$-modules with bounded weight multiplicities, Mosc. Math. J., 6 (2006), 119–134.
\bibitem{LLZ}  G. Liu, R. L\"{u}, K. Zhao, Irreducible Witt modules from Weyl modules and $\gl_n$-modules, J. Algebra, 511 (2018), 164-181.
\bibitem{Liu1}  D. Liu, Y. Pei, L. Xia, Classification of simple weight modules for the $N=2$ superconformal algebra, arXiv:1904.08578.
\bibitem{LX}  R. L\"{u}, Y. Xue, Bounded weight modules over the Lie superalgebra of the Cartan W-type, Algebr. Represent. Theory, (2022),  https://doi.org/10.1007/s10468-021-10112-3.
\bibitem{LZ2}  R. L\"{u}, K. Zhao, Classification of irreducible weight modules over higher rank Virasoro algebras, Adv. Math., 201 (2006), no. 2, 630-656.
\bibitem{Ma}  O. Mathieu, Classification of Harish-Chandra modules over the Virasoro Lie algebras, Invent. Math., 107 (1992), 225-234.
\bibitem{PS}  I. Penkov, V. Serganova, Weight representations of the polynomial Cartan type Lie algebras $W\sb n$ and $\overline S\sb n$, Math. Res. Lett., 6 (1999), no. 3-4, 97-416.
\bibitem{S}  V. Serganova, On representations of Cartan type Lie superalgebras, Translations of the American Mathematical Society-Series 2, 213, (2005), 223-240.
\bibitem{Sh}  G. Shen, Graded modules of graded Lie algebras of Cartan type (I)-Mixed products of modules, Sci. Sinica Ser. A, 29 (1986), no. 6, 570-581.
\bibitem{Su1}  Y. Su, Simple modules over the high rank Virasoro algebras, Commun. Alg., 29 (2001), 2067-2080.
\bibitem{Su2}  Y. Su, Classification of indecomposable $\sl_2(\C)$ modules and a conjecture of Kac on irreducible modules over the Virasoro algebra, J. Algebra, 161 (1993), 33-46.
\bibitem{TZ}  H. Tan, K. Zhao, Irreducible modules over Witt algebras $W_n$ and over $\sl_{n+1}(\C)$, Algebr. Represent. Theory, 21 (2018), no. 4, 787-806.
\bibitem{XL1}  Y. Xue, R. L\"{u}, Classification of simple bounded weight modules of the Lie algebra of vector fields on $\C^n$, arXiv:2001.04204.
\bibitem{XL2}  Y. Xue, R. L\"{u}, Simple weight modules with finite-dimensional weight spaces over Witt superalgebras, J. Algebra 574 (2021), 92–116.
\end{thebibliography}
\end{document}